\documentclass{article}

\usepackage{graphicx,float,latexsym,amssymb,subeqnarray,colordvi}
\usepackage[latin1]{inputenc}
\usepackage{amsmath,amsfonts,amssymb,subeqnarray}

\newtheorem{theorem}{Theorem}

\newtheorem{lemma}[theorem]{Lemma}

\newtheorem{remark}[theorem]{Remark}

\def\R{\mathbb{R}}
\def\C{\mathbb{C}}
\def\dh{{\widehat \delta}}
\def\eh{{\widehat \varepsilon}}

\title{Stability and convergence analysis of
discretizations of the Black--Scholes PDE
\mbox{with the linear boundary condition}}
\author{K.~J.~in 't Hout\footnote{Department of Mathematics and
Computer Science,
University of Antwerp, Middelheimlaan 1, B-2020 Antwerp, Belgium.
\mbox{Email}: \texttt{\{karel.inthout,kim.volders\}@ua.ac.be}.}
~and K. Volders\footnotemark[\value{footnote}]
}

\date{\today}

\begin{document}
\maketitle
\begin{abstract}
In this paper we consider the stability and convergence of numerical
discretizations of the Black--Scholes partial differential equation
(PDE) when complemented with the popular linear boundary condition.
This condition states that the second derivative of the option value
vanishes when the underlying asset price gets large and is often
applied in the actual numerical solution of PDEs in finance.
To our knowledge, the only theoretical stability result in the
literature up to now pertinent to the linear boundary condition
has been obtained by Windcliff, Forsyth \& Vetzal \cite{WFV04}
who showed that for a common discretization a necessary
eigenvalue condition for stability holds.
In this paper, we shall present sufficient conditions for stability
and convergence when the linear boundary condition is employed.
We deal with finite difference discretizations in the spatial
(asset) variable and a subsequent implicit discretization in time.
As a main result we prove that even though the maximum norm of
$e^{tM}$ ($t\ge 0$) can grow with the dimension of the semidiscrete
matrix $M$, this generally does not impair the convergence behavior
of the numerical discretizations.
Our theoretical results are illustrated by ample numerical experiments.
\end{abstract}

\vfill\eject
%%%%%%%%%%%%%%%%%%%%%%%%%%%%%%%%%%%%%%%%%%%%%%%%%%%%%%%%%%%%%%%%%%%%%%%%%%%%%%%%%%%%%%
\section{Introduction}\label{intro}
\noindent
A popular assumption in the valuation of financial options
via the numerical solution of partial differential equations
(PDEs) is the so-called linear boundary condition,
see, for example, \cite{AP10,TR00,W98,WFV04}.
The linear boundary condition states that the second derivative
of the option value with respect to the underlying asset price
vanishes if the asset price gets large.
This condition represents a linear behavior of the option value
for large
%in function of the underlying asset price
asset prices, which can be seen to hold for a wide variety of
financial options.
In spite of its broad use in practice, only few rigorous
results have been derived in the literature up to now on the
stability and convergence of numerical discretizations if
the linear boundary condition is applied.
As it turns out, in the finite difference (FD) approach a
natural treatment of the linear boundary condition leads to
a {\it downwind}\, discretization of the advection term at
the relevant grid point; the details of which are given below
% new
in this section.
Consequently, in the actual numerical solution one might
expect instability, or at least an adverse effect on the
convergence behavior.
% new
It appears, however, that this is not observed in practice.
To our knowledge, the only theoretical stability
analysis in the literature up to now pertinent to the linear
boundary condition has been performed by Windcliff, Forsyth
\& Vetzal \cite{WFV04}.
These authors proved that for a common discretization of the
Black--Scholes PDE a {\it necessary}\, eigenvalue condition
for stability holds.
The objective of the present paper is to arrive at useful
{\it sufficient} conditions for stability and convergence
of discretizations when the linear boundary condition is
employed.
As far as we are aware, such conditions are lacking in the
current literature, but
% new
they are clearly of much interest.
%We consider here the linear boundary condition also in the
%context of the Black--Scholes PDE, so as to explain the
%ideas.
%It is however expected that these ideas extend to more
%advanced, multidimensional PDEs in finance.

Consider the Black--Scholes PDE
\begin{equation}\label{BS}
\frac{\partial u}{\partial t}(s,t) =
\tfrac{1}{2}\sigma^2s^2\frac{\partial^2 u}{\partial s^2} (s,t)
+rs\frac{\partial u}{\partial s}(s,t) - ru(s,t)
\quad (s>0,~ 0< t\leq T),
\end{equation}
where $r>0$ and $\sigma\ge 0$ are given real constants that
denote the risk-neutral interest rate and the volatility,
respectively, and $T>0$ is the given maturity time of the
option.
The exact solution $u(s,t)$ represents the
%new
fair value of an option
if the underlying asset price equals $s$ at time $T-t$.

% new
For the numerical solution, one restricts in practice the
$s$-domain to a boun-ded set
% new
$[0,S]$
with fixed $S>0$ chosen sufficiently large.
The PDE (\ref{BS}) is complemented with initial and
boundary conditions.
% new
In this paper, we consider at $s=S$
the {\it linear boundary condition}
\begin{equation}\label{LBC}
u_{ss}(S,t) = 0 \quad (0< t\leq T).
\end{equation}
At the lower boundary $s=0$ a standard Dirichlet condition is
% new
taken,
which depends on the particular option.
The initial condition is given by the payoff of the option.

FD discretization of the initial-boundary value problem
for (\ref{BS}) on a general (non-uniform) grid
$0=s_0<s_1<s_2<\ldots<s_{m+1}<s_{m+2}=S$, with
mesh widths $h_j=s_j-s_{j-1}$, leads to an initial value
problem for a system of ordinary differential equations
(ODEs),
\begin{equation}\label{ODE}
U'(t)=MU(t)+b(t)\quad (0< t\leq T),\quad U(0)=U_0\,.
\end{equation}
Here $M$ denotes a given real $(m+2)\times(m+2)$--matrix
and $U_0$ and $b(t)$, for $0\leq t\leq T$, are given real
$(m+2)$--vectors.
The vector $U_0$ is
% new
directly given by
the payoff function
and $b$ depends on the Dirichlet condition at $s=0$.
In this paper we
% new
shall
deal with matrices $M$ of the form

\begin{eqnarray}\label{general M}
M
&=&
\left(\begin{array}{c|c}
A&B\\
\hline
&C
\end{array}\right)
\nonumber\\
\nonumber\\
&=&
\left(
\begin{array}{ccccc|cc}
\alpha_1&\gamma_1&&&&&\\
\beta_2&\alpha_2&\gamma_2&&&&\\
&\ddots&\ddots&\ddots&&&\\
&&\beta_{m-1}&\alpha_{m-1}&\gamma_{m-1}&&\\
&&&\beta_{m}&\alpha_{m}&\gamma_{m}&\\
\hline
&&&&&-\frac{rS\phantom{\hat T}}{h_{m+2}}&\frac{rs_{m+1}}{h_{m+2}}\\
&&&&&-\frac{rS\phantom{\hat T}}{h_{m+2}}&\frac{rs_{m+1}}{h_{m+2}}
\end{array}\right)
\end{eqnarray}

\noindent
% new
where $\alpha_j$, $\beta_j$, $\gamma_j$ denote given real numbers.
The
% new
$2\times 2$--matrix $C$
represents a natural discretization of the linear boundary
condition
% new
(\ref{LBC}).
% new
It is determined
by the following approximations at the
grid points $s_{m+1}$ and $s_{m+2}=S$:

\begin{subeqnarray}\label{LBCdiscr}
u_s(s_{m+1},t) \approx \frac{u(s_{m+2},t)-u(s_{m+1},t)}{h_{m+2}}\,,
&& u_{ss}(s_{m+1},t) \approx 0\,, \\
u_s(s_{m+2},t) \approx \frac{u(s_{m+2},t)-u(s_{m+1},t)}{h_{m+2}}\,,
&& u_{ss}(s_{m+2},t) = 0\,.
\end{subeqnarray}
%new
As $r$ in (\ref{BS}) is positive, the
% new
approximation of $u_s$
forms an upwind scheme at $s_{m+1}$, but the same
approximation constitutes a downwind scheme at $s_{m+2}$.
The latter
% new
approximation
can be regarded as obtained from the second-order
central scheme for advection at $s_{m+2}$ with virtual point
$s_{m+3}:=s_{m+2}+h_{m+2}$ and then replacing $u(s_{m+3},t)$
by $2u(s_{m+2},t)-u(s_{m+1},t)$ in view of the linear
boundary condition.
The discretization (\ref{LBCdiscr}b) at $s_{m+2}$
is identical to the one considered in~\cite{WFV04}.
% new
The discretization (\ref{LBCdiscr}a) at $s_{m+1}$, on the
other hand,
appears to be new.
In particular, we approximate $u_{ss}$ by zero at this point
instead of using the
% new
standard second-order central scheme for diffusion.
% new
%This choice has been motivated so as to render our analysis
%feasible.
% new
The choice~(\ref{LBCdiscr}) yields a partial
decoupling between the FD solution at the grid points
$s_{m+1}, s_{m+2}$ and that at $s_1, s_2,\ldots, s_m$.
% new
Concerning
the discretization on $[s_1, s_m]$ we make no assumptions
% new
yet,
except that at each relevant grid point $s_j$ the stencil
belongs to $\{s_{j-1}, s_j, s_{j+1}\}$ -- hence the structure
of the matrices $A$ and $B$ in (\ref{general M}).

% new
The discretization~(\ref{LBCdiscr}) of the linear
boundary condition~(\ref{LBC})
% new
might be interpreted
% new
%and implemented
as a Dirichlet-type condition, since the relevant subsystem
of ODEs involving the matrix $C$ is easily solved exactly
%, given the initial payoff function
% new
(cf. (\ref{expCform}) below).
We
% new
emphasize,
however, that the objective of this paper concerns
comparing the numerical solution to the exact solution of
the Black--Scholes PDE with the linear boundary condition
(\ref{LBC}) at the upper boundary \mbox{$s=S$}, and
{\it not}\, with a Dirichlet condition.
As it turns out, the analysis in the present situation, of
(\ref{LBC}), encounters a variety of additional difficulties.

Our analysis commences with an investigation of the stability
of the FD discretization (\ref{ODE}), (\ref{general M}).
This pertains to the derivation of rigorous bounds on the
norm of the matrix exponential $e^{tM}$.
We shall deal here with the maximum norm.
By $|\cdot|_\infty$ and $\|\cdot\|_\infty$ we denote the
maximum norm of real vectors and matrices, respectively.
An important tool is the {\it logarithmic maximum norm},
which is defined for any square matrix $X$ by
\[
\mu_\infty[X] = \lim_{t\downarrow 0}\frac{||I+tX||_\infty-1}{t}\,,
\]
where $I$ is the identity matrix of the same size as $X$.
Upon writing $X=(\xi_{ij})_{i,j=1}^m$ a convenient formula
for the logarithmic maximum norm is
\begin{equation}\label{mu inf}
\mu_\infty[X]=\max_{1\leq i\leq m} \{ \, \xi_{ii}+\sum_{j\neq i}|\xi_{ij}| \, \}.
\end{equation}
A key property is given by the following theorem;
see~e.g.~\cite{HNW08,HV03,S98,TE05}.

\begin{theorem}\label{stab lognorm}
Let $\omega\in\mathbb{R}$. Then:~~
$
\mu_\infty[X]\leq\omega  ~~\Longleftrightarrow~~
||e^{tX}||_\infty\leq e^{t\omega}~(t\geq 0).
$
\end{theorem}

\noindent
% new
We note that we previously used the logarithmic norm in
analyzing the stability of discretizations of the
Black--Scholes and Heston PDEs when provided with
Dirichlet boundary conditions, see \cite{iHV09,iHV12,V10}.

\vskip0.3cm

An outline of the rest of the paper is as follows.

In Section \ref{stabthm} we investigate the stability of
general semidiscretizations (\ref{ODE}), (\ref{general M})
of the Black--Scholes PDE with the linear boundary condition.
% new
We
prove sharp upper and lower bounds for $||e^{tM}||_\infty$.

In Section \ref{FD} various well-known FD discretizations are
considered. For each discretization a practical sufficient
condition
%on $r$, $\sigma$ and the grid
is
% new
obtained
such that the stability result of Section~\ref{stabthm} holds.

In Section \ref{conv} we
% new
derive a convergence estimate for
general semidiscretizations (\ref{ODE}), (\ref{general M})
of the Black--Scholes PDE with the linear boundary condition.
In Section \ref{stabthm} it was found that $||e^{tM}||_\infty$
% new
is essentially inversely proportional to the mesh width
$h_{m+2}$.
% new
We prove however the positive result that this growth,
as $h_{m+2}$ tends to zero, generally
has no adverse effect on the convergence behavior.

In Section \ref{numexp} extensive numerical experiments
are presented regarding the stability and convergence
results of Sections \ref{stabthm}, \ref{conv}.

%for the five FD schemes that were considered in Section \ref{FD}.

In Section \ref{timediscr} we consider the discretization
in time and prove stability and convergence results for
the popular family of $\theta$-methods.
These results can be
% new
regarded as
%useful
analogues of those obtained for the semidiscretization.

In Section \ref{concl} conclusions and issues for future
research are given.

%%%%%%%%%%%%%%%%%%%%%%%%%%%%%%%%%%%%%%%%%%%%%%%%%%%%%%%%%%%%%%%%%%%%%%%%%%%%%%%%%%%%
\section{A general stability theorem}\label{stabthm}
\setcounter{equation}{0}
\setcounter{theorem}{0}

In this section we consider general matrices $M$ of the form
(\ref{general M}) and derive a useful inclusion for the maximum
norm of $e^{tM}$ for $t\ge 0$.
We start with three lemmas.

\begin{lemma}\label{lemma expM}
It holds that
\[
e^{tM} =
\left(\begin{array}{c|c}
e^{tA}&\int_0^t e^{(t-\tau)A}Be^{\tau C}d\tau\\
\hline
O&e^{tC^{\phantom X}}
\end{array}\right).
\]
\end{lemma}
\textbf{Proof}~
Consider the system of ODEs
\[
U'(t)=M U(t)
\]
with solution given by
\begin{equation}\label{sol}
U(t)=e^{tM}U(0).
\end{equation}
Let the vector $U(t)$ be splitted into two parts,
\[
U(t)=\left(\begin{array}{c}V(t)\\W(t)\end{array}\right),
\]
where $V(t)$ is an $m$--vector and $W(t)$ is a 2--vector.
In view of (\ref{general M}) one has
\[
\begin{cases}
V'(t)\,=AV(t)+BW(t),\\
W'(t)=CW(t).
\end{cases}
\]
Thus $W(t)=e^{tC}W(0)$ and
\begin{eqnarray*}
&&V'(\tau)-AV(\tau)=Be^{\tau C}W(0)\\
&&\Rightarrow\frac{d}{d\tau}\left(e^{-\tau A}V(\tau)\right)=e^{-\tau A}Be^{\tau C}W(0)\\
&&\Rightarrow e^{-tA}V(t)-V(0)=\int_0^te^{-\tau A}Be^{\tau C}W(0)\,d\tau\\
&&\Rightarrow V(t)=e^{tA}V(0)+\int_0^te^{(t-\tau)A}Be^{\tau C}\,d\tau\,\, W(0).
\end{eqnarray*}
Comparing with (\ref{sol}), the result of the lemma is obtained.
\begin{flushright}
$\Box$
\end{flushright}\vskip0.2cm

%note that upper bounds for the norm of $e^{tA}$ have already
%been proven in \cite{iHV09} for second order central schemes
%on non-uniform grids for advection as well as diffusion.

The next lemma gives the maximum norm of $e^{tC}$.
\begin{lemma}\label{lemma expC}
It holds that
\begin{equation}\label{expC}
||e^{tC}||_\infty = e^{-rt}+\left(1-e^{-rt}\right)\frac{2S}{h_{m+2}}\,.
\end{equation}
\end{lemma}
\textbf{Proof}~
The two eigenvalues of $C$ are $0$ and $-r$ with corresponding eigenvectors
$(s_{m+1} ~\, S)^{\rm T}$ and $(1 ~\, 1)^{\rm T}$.
Thus
\begin{equation*}
e^{tC} = \left(\begin{array}{cc}s_{m+1}&1\\S&1\end{array}\right)
\left(\begin{array}{cc}1&0\\0&e^{-rt}\end{array}\right)
\left(\begin{array}{cc}s_{m+1}&1\\S&1\end{array}\right)^{-1}
\end{equation*}
which gives
\begin{equation}\label{expCform}
e^{tC} = \frac{1}{h_{m+2}}\left(\begin{array}{cc}
Se^{-rt}-s_{m+1}&s_{m+1}(1-e^{-rt})\\S(e^{-rt}-1)&S-s_{m+1}e^{-rt}
\end{array}\right).
\end{equation}
Hence,
\[
||e^{tC}||_\infty=
\frac{
\max\left\{|Se^{-rt}-s_{m+1}|+s_{m+1}(1-e^{-rt})\,,\,
S(1-e^{-rt})+S-s_{m+1}e^{-rt}\right\}}{h_{m+2}}\,.
\]
It is readily seen that
\[
|Se^{-rt}-s_{m+1}|+s_{m+1}(1-e^{-rt})\leq
S(1-e^{-rt})+S-s_{m+1}e^{-rt}.
\]
Therefore,
\[
||e^{tC}||_\infty=
\frac{
S(1-e^{-rt})+S-s_{m+1}e^{-rt}}{h_{m+2}}=
e^{-rt}+\left(1-e^{-rt}\right)\frac{2S}{h_{m+2}}.
\]
\begin{flushright}
$\Box$
\end{flushright}\vskip0.2cm

Lemma \ref{lemma expC} shows that for any given $t, r, S >0$ the
maximum norm of $e^{tC}$ is essentially inversely proportional to
the mesh width $h_{m+2}$.
The growth of $||e^{tC}||_\infty$ as $h_{m+2}$ decreases corresponds
to the fact that at the grid point $s=S$ a downwind scheme is used
for the advection term in the Black--Scholes PDE.
% see Section~\ref{intro}.

%We notice that $||e^{tC}||_\infty = |e^{tC}w|_\infty$ with vector
%$w = (-1 ~\, 1)^{\rm T}$.

Let $e_m$ denote the $m$--dimensional unit vector
$(0,\ldots,0,1)^{\rm T}$.

\begin{lemma}\label{lemma Ainvu}
If $A$ is invertible, $\mu_\infty[A]\leq 0$ and
$\alpha_m+|\beta_m|<0$, then
\[
|A^{-1}e_m|_\infty\leq\frac{1}{-\alpha_{m}-|\beta_{m}|}\,.
\]
\end{lemma}
\textbf{Proof}~ (i) Consider first $\mu_\infty[A]<0$.
Define $v=A^{-1}e_m$, i.e., $Av=e_m$.
Writing $v=(v_1,v_2,\ldots,v_m)^{\rm T}$ this gives
the system of equations
\[
\begin{cases}
\alpha_1v_1+\gamma_1v_2=0,&\\
\beta_iv_{i-1}+\alpha_iv_i+\gamma_iv_{i+1}=0&(2\leq i\leq m-1),\\
\beta_{m}v_{m-1}+\alpha_{m}v_{m}=1.&
\end{cases}
\]
We prove by induction that $|v_1|\le |v_2|\le \ldots \le |v_m|$.
In view of (\ref{mu inf}), the assumptions on $A$ imply that
all $\alpha_i <0$ and $\alpha_1+|\gamma_1|<0$.
Using this, yields
\[
|v_1|=\left|\frac{\gamma_1}{\alpha_1}\right|\cdot|v_2|\leq|v_2|.
\]
Next suppose $|v_{i-1}|\leq|v_i|$ for some $2\le i \le m-1$.
We show that $|v_{i}|\leq|v_{i+1}|$.
By (\ref{mu inf}) there holds $|\beta_i|+\alpha_i+|\gamma_i|<0$
and consequently
\[
|\alpha_i|-|\beta_i| >0 \quad{\rm and}\quad
\frac{|\gamma_i|}{|\alpha_i|-|\beta_i|}<1.
\]
We have
\begin{align*}
&\beta_iv_{i-1}+\alpha_iv_i+\gamma_iv_{i+1}=0\\
\Rightarrow\quad
&|\alpha_i|\cdot|v_i|\leq |\beta_i|\cdot|v_{i-1}|
 +|\gamma_i|\cdot|v_{i+1}|\\
\Rightarrow\quad
&|\alpha_i|\cdot|v_i|\leq |\beta_i|\cdot|v_{i}|
 +|\gamma_i|\cdot|v_{i+1}|\\
\Rightarrow\quad
& |v_i|\leq\frac{|\gamma_i|}{|\alpha_i|-|\beta_i|}|v_{i+1}|\\
\Rightarrow\quad
& |v_i|
  \leq|v_{i+1}|.
\end{align*}
This proves the induction step, and it follows that
$|A^{-1}e_m|_\infty = |v|_\infty = |v_{m}|$.
Subsequently,
\begin{align*}
&\alpha_{m} v_{m}=1-\beta_{m} v_{m-1}\\
\Rightarrow\quad
&|\alpha_{m}|\cdot|v_{m}|\leq1+|\beta_{m}|\cdot|v_{m-1}|\\
\Rightarrow\quad
&|\alpha_{m}|\cdot|v_{m}|\leq 1+|\beta_{m}|\cdot|v_{m}|\\
\Rightarrow\quad
&|v_{m}|\leq\frac{1}{|\alpha_{m}|-|\beta_{m}|}\,,
\end{align*}
where it is used that $|\alpha_m|-|\beta_m| >0$.
Since $|\alpha_{m}|=-\alpha_{m}$ the bound of the lemma
is obtained.

(ii) Consider next $\mu_\infty[A]= 0$.
Define the matrix $A_\varepsilon=A-\varepsilon I$ with $\varepsilon>0$.
There holds $\mu_\infty[A_\varepsilon]=-\varepsilon<0$ and thus we can
apply the result from part (i):
\[
|A_\varepsilon^{-1}e_m|_\infty
\leq\frac{1}{-\alpha_{m}+\varepsilon-|\beta_{m}|}\,.
\]
Taking the limit $\varepsilon\downarrow 0$ in this inequality
completes the proof.
\begin{flushright}
$\Box$
\end{flushright}\vskip0.2cm

The following theorem is the first main result of this paper.
It provides a tight inclusion of the maximum norm of $e^{tM}$
for matrices $M$ of the form (\ref{general M}) and reveals
that this norm is essentially inversely proportional to~$h_{m+2}$.

\begin{theorem}\label{general thm}
%Let $M$ be given by {\rm (\ref{general M})}.
If
\begin{equation}\label{condition}
rI+A {\rm~is~invertible},~~\mu_\infty[rI+A]\leq 0,~~
r+\alpha_{m}+|\beta_{m}|+|\gamma_{m}|\leq 0,
\end{equation}
then
\begin{equation}\label{bound}
e^{-rt}+\left(1-e^{-rt}\right)\frac{2S}{h_{m+2}}
\, \leq \, ||e^{tM}||_\infty \, \leq \,
e^{-rt}+\left(1+3e^{-rt}\right)\frac{2S}{h_{m+2}}\,.
\end{equation}
\end{theorem}
\textbf{Proof}~
We employ the formula for $e^{tM}$ given by Lemma \ref{lemma expM}.
First, it is clear that $||e^{tM}||_\infty \ge ||e^{tC}||_\infty$
and by Lemma \ref{lemma expC} the lower bound is directly obtained.
To prove the stated upper bound, we consider the maximum norm of
\[
\int_0^t e^{(t-\tau)A}Be^{\tau C}d\tau.
\]
Using formula (\ref{expCform}), it is readily shown that
\begin{equation}\label{BexpC}
Be^{\tau C} =
\frac{\gamma_{m}}{h_{m+2}}
\left[\begin{array}{cc}(Se^{-r\tau}-s_{m+1})e_m&
(s_{m+1}-s_{m+1}e^{-r\tau})e_m\end{array}\right].
\end{equation}
For any given real numbers  $\phi_0$, $\phi_1$ consider the
vector
\[
f(\tau) = (\phi_0+\phi_1 e^{-r\tau})e_m.
\]
A straightforward computation yields
\[
\int_0^te^{(t-\tau)A}f(\tau)d\tau
=\phi_0 (e^{tA}-I)A^{-1}e_m+
\phi_1 (e^{tA}-e^{-rt}I)(rI+A)^{-1}e_m.
\]
Note that $\mu_\infty[rI+A]\leq 0$ means $\mu_\infty[A]\leq -r$.
In view of this and the assumptions of the theorem it holds
that both $A$ and $rI+A$ are invertible and, by Theorem
\ref{stab lognorm},
\begin{equation}\label{expA}
||e^{tA}||_\infty\leq e^{-rt}.
\end{equation}
Consequently, we have the bound
\[
|\int_0^te^{(t-\tau)A}f(\tau)d\tau\, |_\infty
\leq
(1+e^{-rt})\cdot |\phi_0|\cdot|A^{-1}e_m|_\infty
+2e^{-rt}\cdot|\phi_1|\cdot|(rI+A)^{-1}e_m|_\infty.
\]
Now let $f(\tau)$ represent any column vector of
$Be^{\tau C}$. It is clear from (\ref{BexpC}) that
\[
|\phi_0|\leq |\gamma_{m}|\cdot\frac{S}{h_{m+2}}
\quad {\rm and} \quad
|\phi_1|\leq |\gamma_{m}|\cdot\frac{S}{h_{m+2}}\,.
\]
If $\gamma_m=0$, then the result of the theorem is
obvious; in fact $||e^{tM}||_\infty = ||e^{tC}||_\infty$
in this case.
Thus assume $\gamma_m \not=0$.
Then $r+\alpha_{m}+|\beta_{m}| < 0$ and application of Lemma
\ref{lemma Ainvu} to both $A$ and $rI+A$ yields
\[
|A^{-1}e_m|_\infty\leq\frac{1}{-\alpha_{m}-|\beta_{m}|}
\quad {\rm and} \quad
|(rI+A)^{-1}e_m|_\infty\leq\frac{1}{-r-\alpha_{m}-|\beta_{m}|}\,.
\]
Since
\[
\frac{|\gamma_{m}|}{-\alpha_{m}-|\beta_{m}|} \le 1
\quad {\rm and} \quad
\frac{|\gamma_{m}|}{-r-\alpha_{m}-|\beta_{m}|} \le 1.
\]
it follows that
\[
|\int_0^te^{(t-\tau)A}f(\tau)d\tau\, |_\infty
\leq
(1+3e^{-rt})\cdot \frac{S}{h_{m+2}}\,.
\]
Using that $Be^{\tau C}$ has two columns, we arrive at the bound
\begin{equation}\label{BexpC2}
\| \int_0^t e^{(t-\tau)A}Be^{\tau C}d\tau\, \|_\infty
\le (1+3e^{-rt})\cdot \frac{2S}{h_{m+2}}\,.
\end{equation}
Finally, in view of Lemma \ref{lemma expM},
\[
||e^{tM}||_\infty \le
\max\left\{||e^{tA}||_\infty + \| \int_0^t e^{(t-\tau)A}Be^{\tau C}d\tau\, \|_\infty ~,~
||e^{tC}||_\infty\right\}.
\]
By invoking (\ref{expA}), (\ref{BexpC2}) and Lemma \ref{lemma expC}, the
upper bound of the theorem is obtained.
\begin{flushright}
$\Box$
\end{flushright}\vskip0.2cm

% new
In Section~\ref{FD},
applications of Theorem~\ref{general thm} to
various actual FD discretizations of the Black--Scholes PDE with
the linear boundary condition shall be discussed.
In Section~\ref{conv} the stability results from the present
section shall
% new
effectively be used
in the convergence analysis of FD discretizations.

\begin{remark}
A more direct way to arrive at a bound for $||e^{tM}||_\infty$ is
by using the inequality
\[
\|\int_0^te^{(t-\tau)A}Be^{\tau C}d\tau\, \|_\infty
\leq \int_0^t||e^{(t-\tau)A}||_\infty \cdot ||B||_\infty \cdot ||e^{\tau C}||_\infty\, d\tau\\
%&\leq& |\gamma_m| \int_0^t e^{-r(t-\tau)} \left( e^{-rt}+\left(1-e^{-rt}\right)\frac{2S}{h_{m+2}}\right)
\]
and then applying (\ref{expC}), (\ref{expA}) and calculating the
obtained (simple) integral.
However, this leads to an upper bound which is substantially less
favorable than that of Theorem~\ref{general thm}.
The reason
% new
for this lies in the fact that it yields a factor
$|\gamma_m|/h_{m+2}$ and in all actual applications the quantity
$|\gamma_m|$ is itself inversely proportional to one or more mesh
widths, cf. Section~\ref{FD}.
%We therefore considered the more elaborate construction in the proof
%above, which leads to a sharp upper bound for $||e^{tM}||_\infty$.
%and is key to the positive convergence result of Section \ref{conv}
%below.
\end{remark}

In subsequent applications of Theorem~\ref{general thm} the
following lemma is useful.
\begin{lemma}\label{inverteerbaarheid}
%Let $m\ge 3$.
If $A$ satisfies the conditions
\[
\left\{
  \begin{array}{ll}
    \beta_j \ge 0 & (2\leq j\leq m), \\
    \gamma_j>0 & (1\leq j\leq m-1),\\
    \alpha_1+\gamma_1 \le 0, & \hbox{} \\
    \alpha_j+\beta_j+\gamma_j=0 & (2\leq j\leq m-1), \\
    \alpha_m+\beta_m<0, & \hbox{}
  \end{array}
\right.
\]
then $A$ is invertible and $\mu_\infty[A]=0$.
\end{lemma}
\textbf{Proof}~
The conditions of the lemma directly
% new
imply, by (\ref{mu inf}), that $\mu_\infty[A]=0$.
We
% new
next show that $A$ is invertible.
% new
Upon setting $\beta_1=-\alpha_1-\gamma_1 \ge 0$ and
$\gamma_m=-\alpha_m-\beta_m >0$
%new
and
using that $\alpha_j+\beta_j+\gamma_j\equiv 0$ we obtain
\begin{align*}
&\left(
\begin{array}{cccccc}
\alpha_1&\gamma_1\\
\beta_2&\alpha_2&\gamma_2\\
%&\beta_3&\alpha_3&\gamma_3\\
&\ddots&\ddots&\ddots\\
&&\ddots&\ddots&\ddots\\
&&&\beta_{m-1}&\alpha_{m-1}&\gamma_{m-1}\\
&&&&\beta_m&\alpha_m
\end{array}\right)
\cdot
\left(
\begin{array}{cccccc}
1\\
1&1\\
\vdots&\vdots&\ddots\\
\vdots&\vdots&&\ddots\\
%1&&&&\ddots\\
1&1&\ldots&\ldots&1
\end{array}
\right)=\\
&\left(
\begin{array}{cccccc}
-\beta_1&\gamma_1\\
&-\beta_2&\gamma_2\\
&&\ddots&\ddots\\
&&&\ddots&\ddots\\
&&&&-\beta_{m-1}&\gamma_{m-1}\\
\hline
-\gamma_m&-\gamma_m&\ldots&\ldots&-\gamma_m&\alpha_m
\end{array}
\right)=:P\,.
\end{align*}
Next, modify the matrix $P$ by subtracting column 1 from
columns $2, 3, \ldots, m-1$.
This leads to
\[\widetilde{P}=
\left(
\begin{array}{cccccc|c}
-\beta_1&-\alpha_1&\beta_1&\beta_1&\ldots&\beta_1&\\
\hline
&-\beta_2&\gamma_2&&&&\\
&&-\beta_3&\gamma_3&&&\\
&&&\ddots&\ddots&&\\
&&&&-\beta_{m-2}&\gamma_{m-2}&\\
&&&&&-\beta_{m-1}&\gamma_{m-1}\\
\hline
-\gamma_m&&&&&&\alpha_m
\end{array}
\right).
\]
Clearly, if $\widetilde{P}$ is invertible, then so is $A$.
We prove that
\[
x\in\mathbb{R}^m\,,\, \widetilde{P}x=0
~\Longrightarrow~ x=0.
\]
Write $x=(x_1,x_2,\ldots,x_m)^{\rm T}$. Then $\widetilde{P}x=0$
is equivalent to the system of equations
\begin{equation}\label{vgln}
\begin{cases}
-\beta_1x_1-\alpha_1x_2+\beta_1x_3+\beta_1x_4+\ldots+\beta_1x_{m-1}\phantom{+\gamma_{m-1}x_m}
&=\, 0\\
\phantom{-\beta_1x_1}-\beta_2x_2+\gamma_2x_3\phantom{+\beta_1x_4+\ldots+\beta_1x_{m-1}}
&=\, 0\\
\phantom{-\beta_1x_1-\alpha_1x_2}-\beta_3x_3+\gamma_3x_4\phantom{+\ldots+\beta_1x_{m-1}}
&=\, 0\\
\phantom{-\beta_1x_1-\alpha_1)_2-\beta_3x_3+\gamma_3x_4}\quad\ddots
&\quad\vdots\\
\phantom{-\beta_1x_1-\alpha_1x_2+\beta_1x_3+\beta_1x_4+\ddots}-\beta_{m-1}x_{m-1}+\gamma_{m-1}x_m
&=\, 0\\
-\gamma_mx_1\phantom{-\alpha_1x_2+\beta_1x_3+\beta_1x_4+\ldots-\beta_{m-1}x_{m-1}+}+\alpha_mx_m
&=\, 0.
\end{cases}
\end{equation}

\noindent
We distinguish three cases.

(a) Assume $\beta_j \not= 0$ whenever $1\le j\le m$.
Then $\alpha_j<0$, $\beta_j>0$, $\gamma_j>0$ for all~$j$.
Starting with the last equation of (\ref{vgln}) and moving upwards,
one finds that
\begin{equation*}
x_m =
\frac{\gamma_m}{\alpha_m}\,x_1\,,\quad
x_{j} =
\left(\prod_{k=j}^{m-1}\frac{\gamma_k}{\beta_k}\right)\frac{\gamma_m}{\alpha_m}\,x_1
\qquad (2\leq j\leq m-1).
\end{equation*}
Substituting this into the first equation of (\ref{vgln}) yields
\[
\left[
-\beta_1
-\alpha_1\left(
          \prod_{k=2}^{m-1}\frac{\gamma_k}{\beta_k}\right)
          \frac{\gamma_m}{\alpha_m}
+\beta_1\sum_{j=3}^{m-1}\left(
                         \prod_{k=j}^{m-1}\frac{\gamma_k}{\beta_k}
                         \right)
                         \frac{\gamma_m}{\alpha_m}
\right]
x_1
=0.
\]
It is easily seen that the coefficient of $x_1$ in the latter equation
is nonzero.
Thus $x_1=0$, and consequently $x=0$.

(b) Assume $\beta_1 = 0$. Since $\alpha_1$ and all $\gamma_j$
are nonzero, (\ref{vgln}) directly yields that $x=0$.

(c) Assume $\beta_j = 0$ for certain $2\le j\le m$.
This induces a natural partitioning of the matrix $A$ where each
diagonal block belongs to either case (a) or case (b) above.
Using this, it readily follows that if $Ax=0$ then $x=0$.
\begin{flushright}
$\Box$
\end{flushright}\vskip0.2cm

\begin{remark}
Under the
% new
assumptions
of Lemma~\ref{inverteerbaarheid} it holds
that $-A$ is a so-called
% new
M-matrix.
% new
This follows, e.g.,
by using condition $({\rm M}_{37})$ in
\cite[Chapter~6]{BP94}.
\end{remark}

%%%%%%%%%%%%%%%%%%%%%%%%%%%%%%%%%%%%%%%%%%%%%%%%%%%%%%%%%%%%%%%%%%%%%%%%%%%%%%%%%%%%%%
\section{Application of the general stability analysis}\label{FD}
\setcounter{equation}{0}
\setcounter{theorem}{0}

In the following we apply the general stability analysis of Section
\ref{stabthm} to actual FD discretizations of the Black--Scholes PDE
with the linear boundary condition.
% new
%For each semidiscretization a practical condition on $r$, $\sigma$
%and the grid is obtained such that the stability bound (\ref{bound})
%holds.
Pertinent to the advection term
% new
on $[s_1,s_m]$ we consider three FD schemes:
\begin{subeqnarray}\label{FDadv}
u_s(s_j,t)&\approx&
\frac{u(s_{j+1},t)-u(s_j,t)}{h_{j+1}}\,,\\
u_s(s_j,t)&\approx&\frac{u(s_{j+1},t)-u(s_{j-1},t)}{H_j}\,,\\
u_s(s_j,t)&\approx&-\frac{h_{j+1}}{h_jH_j} u(s_{j-1},t) +
\frac{h_{j+1}-h_j}{h_jh_{j+1}} u(s_j,t)+
\frac{h_j}{h_{j+1}H_j} u(s_{j+1},t).\qquad~~
\end{subeqnarray}
The scheme (\ref{FDadv}a) has a first-order truncation
error and is called the \emph{first-order forward scheme}.
The scheme (\ref{FDadv}b) possesses a second-order truncation
error whenever the grid is smooth; for the scheme (\ref{FDadv}c)
this holds
% new
for arbitrary grids.
We refer to (\ref{FDadv}b) and (\ref{FDadv}c) as the
\emph{central scheme~A} and the \emph{central scheme~B},
respectively.
%(cf.~e.g.~\cite{iHV09,VR92}).
Notice that these two schemes are identical if the grid
is uniform.
For the diffusion term
% new
on $[s_1,s_m]$
we use the standard central FD scheme
\begin{equation}\label{FDdiff}
u_{ss}(s_j,t)\approx
\frac{2}{h_{j}H_j}u(s_{j-1},t)-\frac{2}{h_jh_{j+1}}u(s_j,t)
+\frac{2}{h_{j+1}H_j}u(s_{j+1},t)\,,
\end{equation}
which possesses a second-order truncation error on smooth
grids.

Based on (\ref{FDadv}), (\ref{FDdiff}) we consider five
well-known FD discretizations of the Black--Scholes PDE.
These discretizations differ in the treatment of the
advection term $rsu_s$ on the interval $[s_1,s_m]$.
The linear boundary condition is always discretized on
$[s_{m+1},s_{m+2}]$ according to (\ref{LBCdiscr}).
For each semidiscretization, successive application of
Lemma \ref{inverteerbaarheid} (with $A$ replaced by $rI+A$)
and Theorem~\ref{general thm} directly
% new
gives a
condition on $r$, $\sigma$ and the grid such that the
stability bound (\ref{bound})
% new
holds.

\begin{enumerate}

\item[1.] {\it Forward.}~~Using the first-order
forward scheme gives\footnote{Note that $\beta_1$
need not to be defined, but it is convenient for
the analysis.} for $1\le j\le m$,
\begin{equation}\label{forward + diffusion}
\begin{cases}
\beta_j=
\sigma^2\displaystyle\frac{s_j^2}{h_jH_j}\,,\\
\alpha_j=
-r-r\displaystyle\frac{s_j}{h_{j+1}}-
\sigma^2\displaystyle\frac{s_j^2}{h_jh_{j+1}}\,,\\
\gamma_j=
r\displaystyle\frac{s_j}{h_{j+1}}+\sigma^2\frac{s_j^2}{h_{j+1}H_j}\,.
\end{cases}
\end{equation}
The bound (\ref{bound}) holds for all $r>0$, $\sigma\ge 0$ and
all grids.

\item[2.] {\it Central A.}~~Using the central scheme A
gives for $1\le j\le m$,
\begin{equation}\label{method A + diffusion}
\begin{cases}
\beta_j=-r\displaystyle\frac{s_j}{H_j}+
\sigma^2\displaystyle\frac{s_j^2}{h_jH_j}\,,\\
\alpha_j=-r-\sigma^2\displaystyle\frac{s_j^2}{h_jh_{j+1}}\,,\\
\gamma_j=r\displaystyle\frac{s_j}{H_j}+
\sigma^2\displaystyle\frac{s_j^2}{h_{j+1}H_j}\,.
\end{cases}
\end{equation}
The bound (\ref{bound}) holds if for all $1\le j\le m$:
\begin{equation}\label{sufA}
0<r\le \frac{s_j}{h_j}\sigma^2.
\end{equation}

\item[3.] {\it Central B.}~~Using the central scheme B
gives for $1\le j\le m$,
\begin{equation}\label{method B + diffusion}
\begin{cases}
\beta_j=-r\displaystyle\frac{s_jh_{j+1}}{h_jH_j}+
\sigma^2\displaystyle\frac{s_j^2}{h_jH_j}\,,\\
\alpha_j=-r+r\displaystyle\frac{s_j(h_{j+1}-h_j)}{h_jh_{j+1}}-
\sigma^2\displaystyle\frac{s_j^2}{h_jh_{j+1}}\,,\\
\gamma_j=r\displaystyle\frac{s_jh_j}{h_{j+1}H_j}+
\sigma^2\displaystyle\frac{s_j^2}{h_{j+1}H_j}\,.
\end{cases}
\end{equation}
The bound (\ref{bound}) holds if for all $1\le j\le m$:
\begin{equation}\label{sufB}
0<r\le \frac{s_j}{h_{j+1}}\sigma^2.
\end{equation}

\item[4.] {\it Mixed A.}~~This scheme is defined as a
suitable combination of Forward and Central A above.
For each given $1\le j\le m$: if (\ref{sufA}) is
fulfilled, then define $\alpha_j$, $\beta_j$, $\gamma_j$
by (\ref{method A + diffusion}); else by
(\ref{forward + diffusion}).
This scheme has been considered e.g.~in \cite{WFV04}.
Clearly, the idea is to switch from the central scheme
(\ref{method A + diffusion}) to the forward scheme
(\ref{forward + diffusion}) in those grid points that
would give rise to a strictly positive logarithmic
maximum norm of $rI+A$.
By construction, it follows that the bound (\ref{bound})
holds for all $r>0$, $\sigma\ge 0$ and all grids.

\item[5.] {\it Mixed B.}~~This scheme is defined as a
suitable combination of Forward and Central B above.
For each given $1\le j\le m$: if (\ref{sufB}) is
fulfilled, then define $\alpha_j$, $\beta_j$, $\gamma_j$
by (\ref{method B + diffusion}); else by
(\ref{forward + diffusion}).
Again, it follows that the bound (\ref{bound})
holds for all $r>0$, $\sigma\ge 0$ and all grids.
\end{enumerate}

%\section{Stability estimates in a scaled Euclidean norm}
%\setcounter{equation}{0}
%\setcounter{theorem}{0}
%
%In this section we give rigorous upper bounds for
%$||e^{tM}||_H$. We can make use of the stability
%results for the maximum norm by using the following Lemma.
%
%\begin{lemma} Let $R$ be a given $(m+2)\times(m+2)$--matrix.
%It holds that
%\[
%||R||_H\leq
%\sqrt{\frac{\sum_{i=1}^{m+2}H_i}{\min_{1\leq i\leq m+2}\{H_i\}}}\cdot||R||_\infty.
%\]
%\end{lemma}
%\textbf{Proof}~
%Let $x\in\mathbb{R}^{m+2}$, then
%\[
%|x|_H=\sqrt{\sum_{i=1}^{m+2}\frac{H_i}{2}|x_i|^2}\,
%\leq\left(\sqrt{\sum_{i=1}^{m+2}\frac{H_i}{2}}\,\,\right)\cdot|x|_\infty
%\]
%and
%\[
%|x|_H=\sqrt{\sum_{i=1}^{m+2}\frac{H_i}{2}|x_i|^2}\,
%\geq\left(\sqrt{\frac{\min_{1\leq i\leq m+2}H_i}{2}}\,\,\right)\cdot|x|_\infty
%\]
%Thus
%\begin{eqnarray*}
%||R||_H&=&\max_{x\neq 0}\frac{|Rx|_H}{|x|_H}\\
%&\leq&\max_{x\neq0}\frac{\sqrt{\sum_{i=1}^{m+2}\frac{H_i}{2}}\cdot|Rx|_\infty}
%{\sqrt{\frac{\min_{1\leq i\leq m+2}\{H_i\}}{2}}\cdot|x|_\infty}\\
%&=&\sqrt{\frac{\sum_{i=1}^{m+2}{H_i}}{{\min_{1\leq i\leq m+2}\{H_i\}}}}\cdot|R|_\infty
%\end{eqnarray*}
%\begin{flushright}
%$\Box$
%\end{flushright}\vskip0.2cm

%%%%%%%%%%%%%%%%%%%%%%%%%%%%%%%%%%%%%%%%%%%%%%%%%%%%%%%%%%%%%%%%%%%%%%%%%%%%%%%%%%%%%%
\section{A general convergence result}\label{conv}
\setcounter{equation}{0}
\setcounter{theorem}{0}

In this section we
% employ the stability analysis of Section~\ref{stabthm} to
prove a convergence result for general semidiscretizations
(\ref{ODE}), (\ref{general M}) of the Black--Scholes PDE
with the linear boundary condition.
Let~$U$ be the exact solution to (\ref{ODE}),
(\ref{general M}) and, for $0\le t\le T$, let $u_h(t)$
be the vector of the same size as $U(t)$ given by
\[
u_h(t)=(u(s_1,t),u(s_2,t),\ldots,u(s_{m+2},t))^{\rm T}\,,
\]
where $u$ is the exact solution to the initial-boundary
value problem for the Black--Scholes PDE (\ref{BS}) on
$0\le s\le S$ with linear boundary condition (\ref{LBC}).
Define the \emph{spatial discretization error}
\begin{equation*}\label{global error}
\varepsilon_h(t)=u_h(t)-U(t)
\end{equation*}
and the \emph{spatial truncation error}
\begin{equation*}\label{local error}
\delta_h(t)=u_h'(t)-Mu_h(t)-b(t).
\end{equation*}
% new
A
standard approach to estimate spatial discretization
errors is to combine an estimate for the spatial truncation
errors with a stability bound, cf.~e.g.~\cite{HV03}.
However, a direct use of the bound on $||e^{tM}||_\infty$
from Theorem~\ref{general thm} does not lead to an optimal
result.
To obtain a useful result in the present situation where
the linear boundary condition is employed, we consider a
partitioning of the spatial truncation error vector into
two parts, corresponding to the intervals $[s_1,s_m]$ and
$[s_{m+1},s_{m+2}]$:
\[
\delta_h(t)=
\left(
  \begin{array}{c}
    \delta_h^L(t)\vspace{0.1cm} \\
    \delta_h^R(t) \\
  \end{array}
\right)
\quad
{\rm with}~~\delta_h^L(t)\in\R^m~~{\rm and}~~\delta_h^R(t)\in\R^2.
\]
Using
% new
the
individual stability bounds derived in Section \ref{stabthm},
% new
we have as a preliminary result
\begin{lemma}\label{epsilon}
%Let $M$ be given by {\rm (\ref{general M})} and satisfy condition
Assume {\rm(\ref{condition})} holds.
Then the spatial discretization error satisfies
\[
|\varepsilon_h(t)|_\infty
\leq t\cdot
\max_{0\leq\vartheta\leq t}
\left\{
|\delta_h^L(\vartheta)|_\infty+\frac{8S}{h_{m+2}}\cdot|\delta_h^R(\vartheta)|_\infty
\right\}
\quad (t>0).
\]
\end{lemma}
\textbf{Proof}~
From
\begin{eqnarray*}
u_h'(t)&=&Mu_h(t)+b(t)+\delta_h(t),\\
U'(t)&=&MU(t)+b(t)
\end{eqnarray*}
one has
\[
\varepsilon_h'(t)=M\varepsilon_h(t)+\delta_h(t)~~~~~~~~~
\]
and, since $\varepsilon_h(0)=0$,
\[
\varepsilon_h(t)=\int_0^te^{(t-\vartheta)M}\delta_h(\vartheta)d\vartheta.
\]
Lemma \ref{lemma expM} yields
\[
e^{(t-\vartheta)M}\delta_h(\vartheta)
=
\left(\begin{array}{c|c}
e^{(t-\vartheta)A}&\int_0^{t-\vartheta} e^{(t-\vartheta-\tau)A}Be^{\tau C}d\tau\\
\hline
O&e^{{(t-\vartheta)C}^{\phantom X}}
\end{array}\right)
\left(\begin{array}{c}
\delta_h^L(\vartheta)\vspace{0.1cm}\\
\delta_h^R(\vartheta)
\end{array}
\right).
\]
By (\ref{expA}), (\ref{BexpC2}) and Lemma \ref{lemma expC} the
following bounds hold whenever $0\le \vartheta \le t$\,:
\begin{eqnarray*}
||e^{(t-\vartheta)A}||_\infty&\leq &1,\\
||\int_0^{t-\vartheta}e^{(t-\vartheta-\tau)A}Be^{\tau C}d\tau||_\infty
&\leq &\frac{8S}{h_{m+2}}\,,\\
||e^{(t-\vartheta)C}||_\infty&\leq&\frac{2S}{h_{m+2}}\,.
\end{eqnarray*}
Hence
\[
|e^{(t-\vartheta)M}\delta_h(\vartheta)|_\infty \le
|\delta_h^L(\vartheta)|_\infty+\frac{8S}{h_{m+2}}\cdot|\delta_h^R(\vartheta)|_\infty.
\]
Together with the integral representation above, this
gives the bound on the maximum norm of $\varepsilon_h(t)$.
\begin{flushright}
$\Box$
\end{flushright}\vskip0.2cm

The following theorem forms the second main result of this
paper.
It gives a useful convergence estimate for general
semidiscretizations (\ref{ODE}), (\ref{general M})
of the Black--Scholes PDE with linear boundary condition.
% new
Its
proof is obtained by combining Lemma \ref{epsilon} with
a bound for $\delta_h^R$.

\begin{theorem}\label{convergence}
Let $h^\ast> 0$ be given and assume that on $[S-h^\ast,S]\times [0,T]$
the partial derivative $u_{sss}$ exists and is continuous.
Define
\begin{equation}\label{kappa_eta}
\kappa  = 4\sigma^2 S^3 + 6r S^2 h^\ast~~{\it and}~~
\eta(t) = \max \{ |u_{sss}(\xi,t)|: S-h^\ast\le \xi \le S \}.
\end{equation}
%Let $M$ be given by {\rm (\ref{general M})} and satisfy condition
Assume {\rm(\ref{condition})} holds.
Then the spatial discretization error satisfies
\[
|\varepsilon_h(t)|_\infty
\leq t\cdot
\max_{0\leq\vartheta\leq t}
\left\{
|\delta_h^L(\vartheta)|_\infty+\kappa\eta(\vartheta)
\right\}
\quad {\it whenever}~ 0<h_{m+2} \le h^\ast,~ 0<t\le T.
\]
\end{theorem}
\textbf{Proof}~Write $s=s_{m+1}$ and $h=h_{m+2}$.
Pertinent to the point $s_{m+1}$ we have
\begin{eqnarray*}
\delta_{h,1}^R(t)
&=& u_t(s,t) - r s\, \frac{u(S,t)-u(s,t)}{h} + r u(s,t) \\
&=& \tfrac{1}{2}\sigma^2 s^2 u_{ss}(s,t)+
r s \, \left[ u_s(s,t) - \frac{u(S,t)-u(s,t)}{h} \right].
\end{eqnarray*}
By Taylor's theorem and using $u_{ss}(S,t)=0$ there
follows
\begin{eqnarray*}
u(s,t)      &=& u(S,t) - hu_s(S,t) + E_0(h,t),  \\
u_s(s,t)    &=& u_s(S,t) + E_1(h,t), \\
u_{ss}(s,t) &=& E_2(h,t),
\end{eqnarray*}
where
\[
E_0(h,t) = -\tfrac{1}{6} h^3\, u_{sss}(\xi_0,t)~,~~
E_1(h,t) = \tfrac{1}{2} h^2\, u_{sss}(\xi_1,t)~,~~
E_2(h,t) = -h\, u_{sss}(\xi_2,t)
\]
with certain $\xi_0, \xi_1, \xi_2$ in $(s,S)$.
Substitution into the above formula yields
\[
\delta_{h,1}^R(t) = \tfrac{1}{2}\sigma^2 s^2 E_2(h,t)
+rs  E_1(h,t) +  \frac{rs}{h} E_0(h,t),
\]
which readily leads to the estimate
\begin{equation*}\label{d1}
|\delta_{h,1}^R(t)| \le \{ \tfrac{1}{2}\sigma^2 S^2
+ \tfrac{2}{3}rS h \}\eta(t) \cdot h.
\end{equation*}
\noindent
Analogously, pertinent to the point $s_{m+2}$, there holds
\begin{eqnarray*}
\delta_{h,2}^R(t)
&=& u_t(S,t) - r S\, \frac{u(S,t)-u(s,t)}{h} + r u(S,t) \\
&=& r S \, \left[ u_s(S,t) - \frac{u(S,t)-u(s,t)}{h} \right]\\
&=& \frac{rS}{h} E_0(h,t)
\end{eqnarray*}
and
\begin{equation*}\label{d2}
|\delta_{h,2}^R(t)| \le \tfrac{1}{6}rS \eta(t)\cdot h^2.
\end{equation*}
It thus follows that
\begin{equation*}
|\delta_h^R(t)|_\infty \le \{ \tfrac{1}{2}\sigma^2 S^2
+ \tfrac{2}{3}rS h^\ast\}\eta(t)\cdot h
\end{equation*}
and
\begin{equation}\label{deltahR}
\frac{8S}{h}\cdot|\delta_h^R(t)|_\infty \le
\kappa \eta(t).
\end{equation}
Application of Lemma \ref{epsilon} then gives the desired
estimate for $|\varepsilon_h(t)|_\infty$.
\begin{flushright}
$\Box$
\end{flushright}\vskip0.2cm

The estimate of Theorem~\ref{convergence} for the spatial
discretization error consists of two contributions,
corresponding to the two intervals $[s_1,s_m]$ and
$[s_{m+1},s_{m+2}]$.
The first contribution is equal to the part of the spatial
truncation error pertinent to $[s_1,s_m]$.
For any given FD discretization this can be estimated in a
standard way by Taylor expansion.
The second contribution is equal to $\kappa \eta(\vartheta)$
and depends on the partial derivative $u_{sss}$ of the exact
option price near $s=S$.
For a wide range of financial options it is plausible that
this contribution can be made arbitrarily small upon taking
the upper bound $S$ sufficiently
% new
large (in the case of European call and put options this is
readily proved).
Theorem~\ref{convergence} thus expresses the useful result
that the contribution to spatial discretization error of the
semidiscretized linear boundary condition is
% new
negligible, provided $S$ is sufficiently large.
Key to the proof is that the stability bounds derived in
Section~\ref{stabthm} admit a growth of errors from the interval
$[s_{m+1},s_{m+2}]$ that is at most inversely proportional to
$h_{m+2}$ but this growth is precisely offset by the factor
$h_{m+2}$ that arises in the part of the spatial truncation
error pertinent to this interval.

We note that
% new
in Theorem~\ref{convergence} it is tacitly assumed
that the initial function is smooth, which has been
used for ease of
% new
the analysis, but is often not fulfilled in
% new
applications.
The numerical experiments in the subsequent section deal
with a nonsmooth initial function.

%is governed by the
%spatial truncation error of the semidiscretization on the
%interval $[s_1,s_m]$ whenever $S$ is sufficiently large.

%From this we can conclude that, although the upper bound for
%$||e^{tM}||_\infty$ in Theorem \ref{stab forward} grows with the
%dimension of the matrix $M$, we still maintain first order convergence.
%The reason for this is that the elements of the matrix $M$ that
%cause this instability (i.e. the last two columns of $M$) correspond
%to the grid points where the linear boundary condition holds, and therefore
%the corresponding spatial truncation error have order 2 instead of 1.
%These two effects cancel each other out.

%%%%%%%%%%%%%%%%%%%%%%%%%%%%%%%%%%%%%%%%%%%%%%%%%%%%%%%%%%%%%%%%%%%%%%%%%%%%%%%%%%%%%%
\section{Numerical experiments}\label{numexp}
\setcounter{equation}{0}
\setcounter{theorem}{0}

In this section the stability and convergence results
of Sections~\ref{stabthm},~\ref{conv} are illustrated
by
% new
%ample
numerical experiments.
We consider the five FD discretizations of the
Black--Scholes PDE formulated in
% new
Section~\ref{FD} with
the linear boundary condition discretized
as given in Section~\ref{intro}.
For each FD discretization, we also consider its
analogue where the definition of $\alpha_j$,
$\beta_j$, $\gamma_j$ for $1\le j\le m$ in Section~\ref{FD}
is extended to $j=m+1$. This corresponds to the approach
by Windcliff, Forsyth \& Vetzal \cite{WFV04} where in the
penultimate grid point $s_{m+1}$ no modified discretization
is applied related to the linear boundary condition.
In the following we shall refer by LBC1 to the numerical
treatment of the linear boundary condition as defined in
Section~\ref{intro} and by LBC2 to the treatment as
considered in \cite{WFV04}.

For the numerical experiments a typical smooth, non-uniform
spatial grid is chosen.
Let $E\in (0,S)$ and $c>0$ be given and fixed.
Consider the continuous, strictly increasing function
\begin{equation*}
\varphi (\xi) = E + c\cdot \sinh(\xi) \quad(a\le \xi \le b)
\end{equation*}
with $a=\sinh^{-1}(-E/c)$ and $b=\sinh^{-1}((S-E)/c)$.
Let
\[
\xi_j = a+j\cdot\Delta \xi \quad (0\leq j\leq m+2)
~~~\textrm{with}~~ \Delta \xi = \frac{b-a}{m+2}.
\]
Then a non-uniform grid $0=s_0 < s_1 < \ldots < s_{m+1} < s_{m+2}=S$
is defined by the transformation
\begin{equation*}
\label{grid} s_j=\varphi(\xi_j) \quad (0\leq j\leq m+2).
\end{equation*}
The parameter $E$ can be viewed as the exercise price of a vanilla
option and $c$ determines the fraction of grid points $s_j$ that
lie in the neighborhood of $E$,
\[
h_j\approx c\cdot\Delta\xi \quad \textrm{whenever}~~s_j\approx E.
\]
The grid is
% new
{\it smooth}
in the sense that there exist real constants
$c_0$, $c_1$, $c_2>0$ (independent of $j$ and $m$) such that the
mesh widths $h_j$ satisfy
\begin{equation*}
\label{smooth} c_0\cdot\Delta\xi\leq h_j\leq c_1\cdot\Delta \xi
\quad \textrm{and} \quad |h_{j+1}-h_j|\leq c_2\cdot(\Delta\xi)^2.
\end{equation*}
The above type of grid is often used in financial applications,
cf.~e.g.~\cite{iHV09,TR00}.
We (arbitrarily) set $E=100$, $c=E/5$, $S=400$.

\subsection{Stability experiments}\label{stabexp}
\begin{figure}
\begin{center}
\begin{tabular}{cc}
\includegraphics[width=0.5\textwidth]{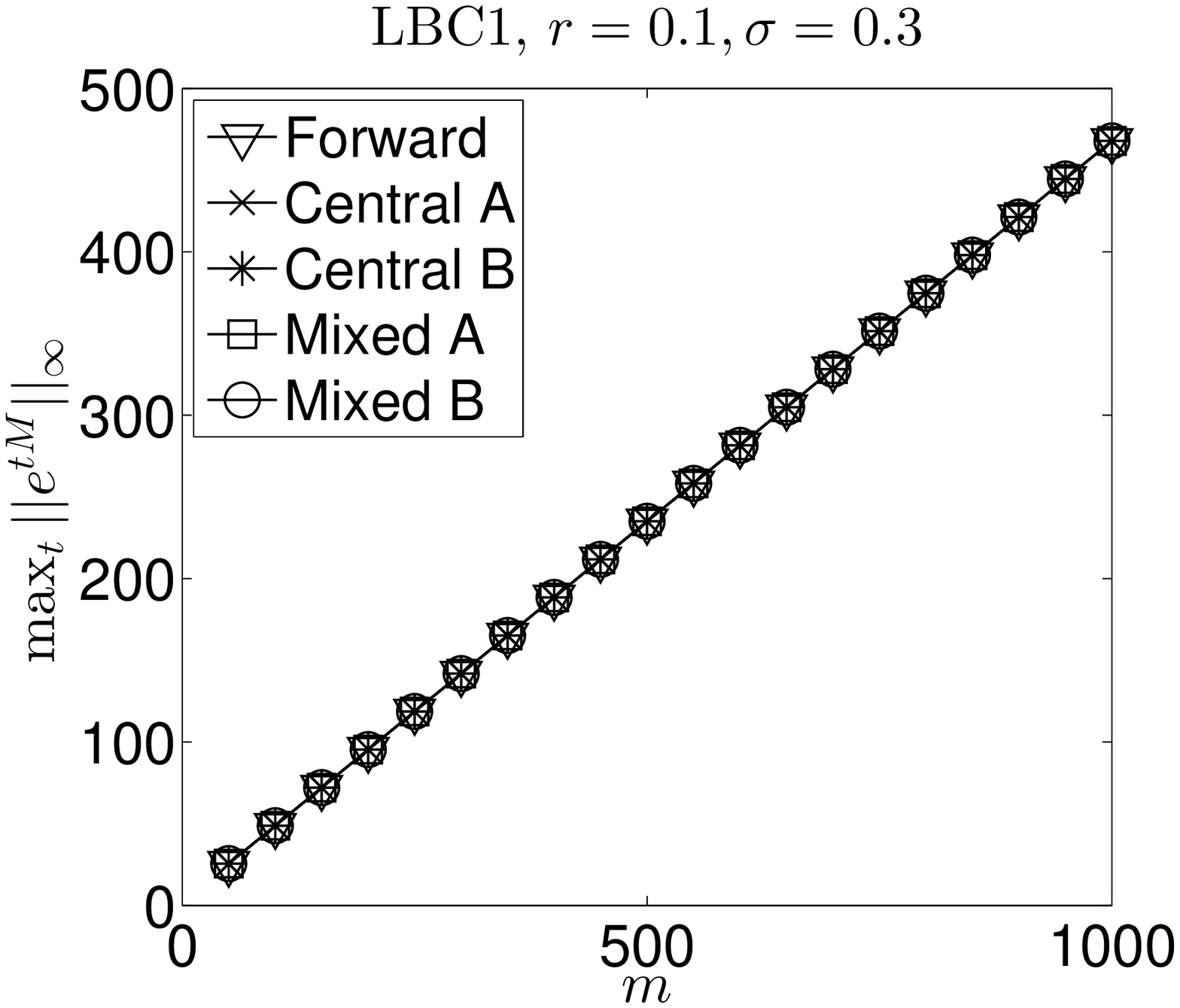}
&\includegraphics[width=0.5\textwidth]{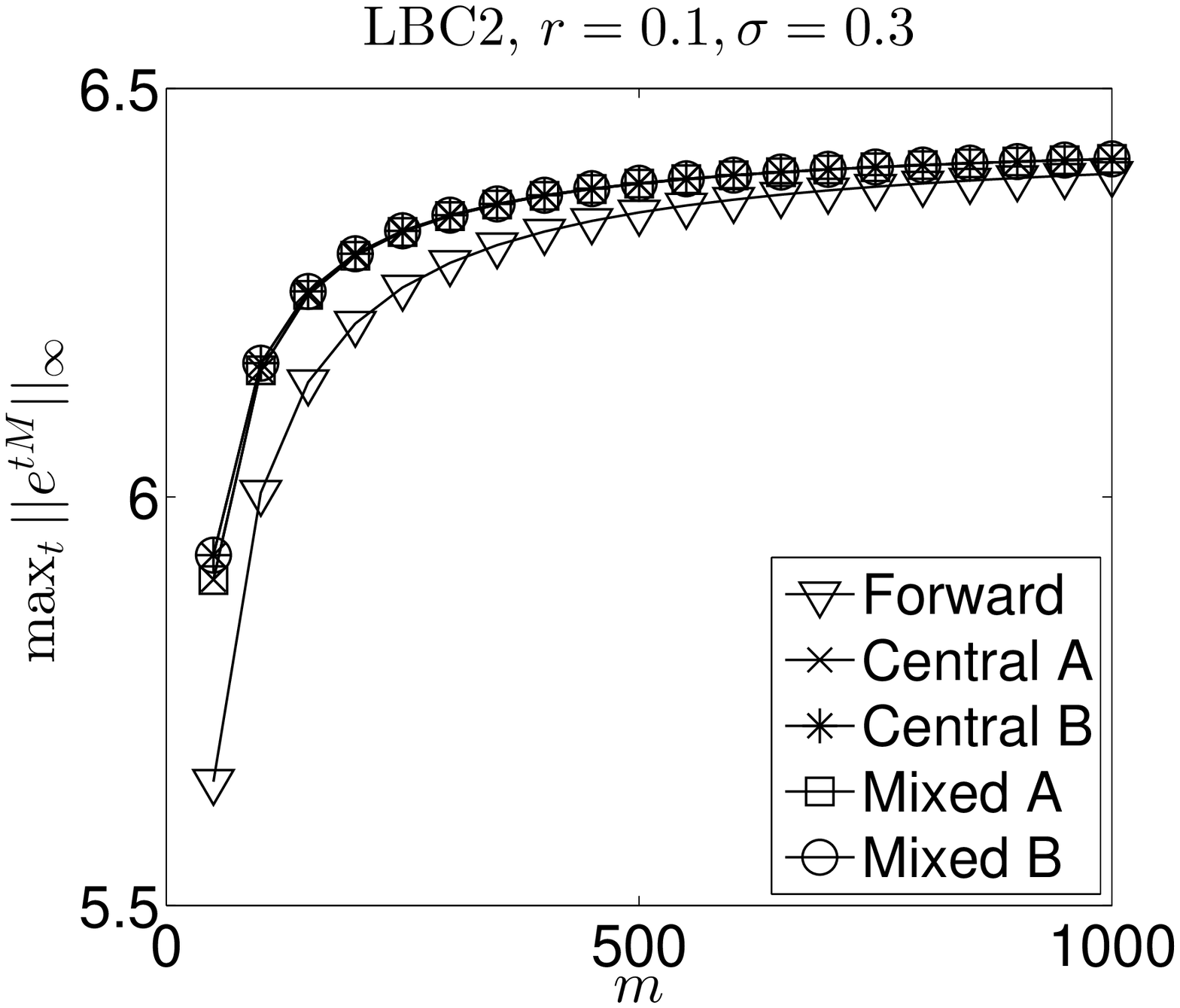}\\
\includegraphics[width=0.5\textwidth]{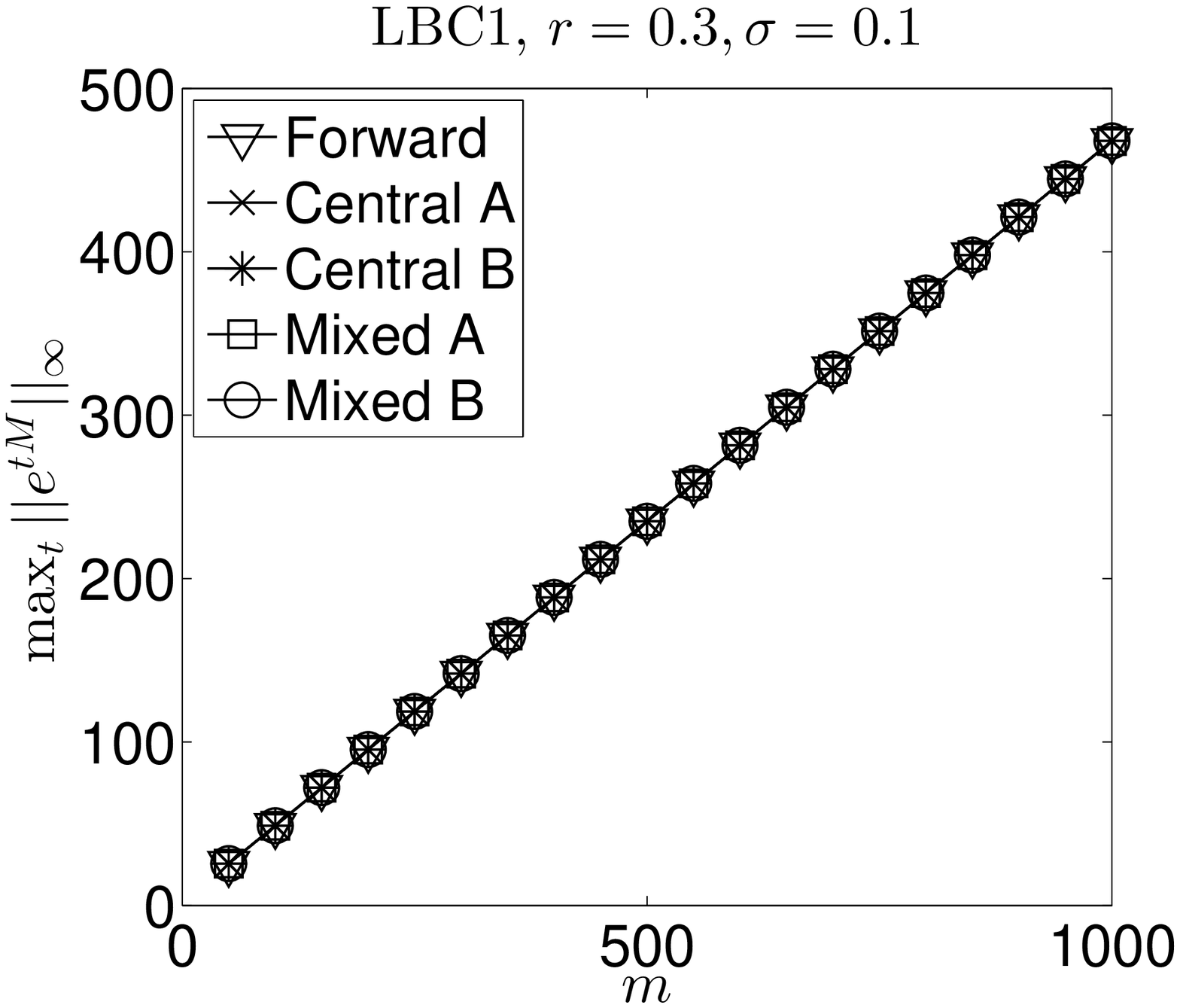}
&\includegraphics[width=0.5\textwidth]{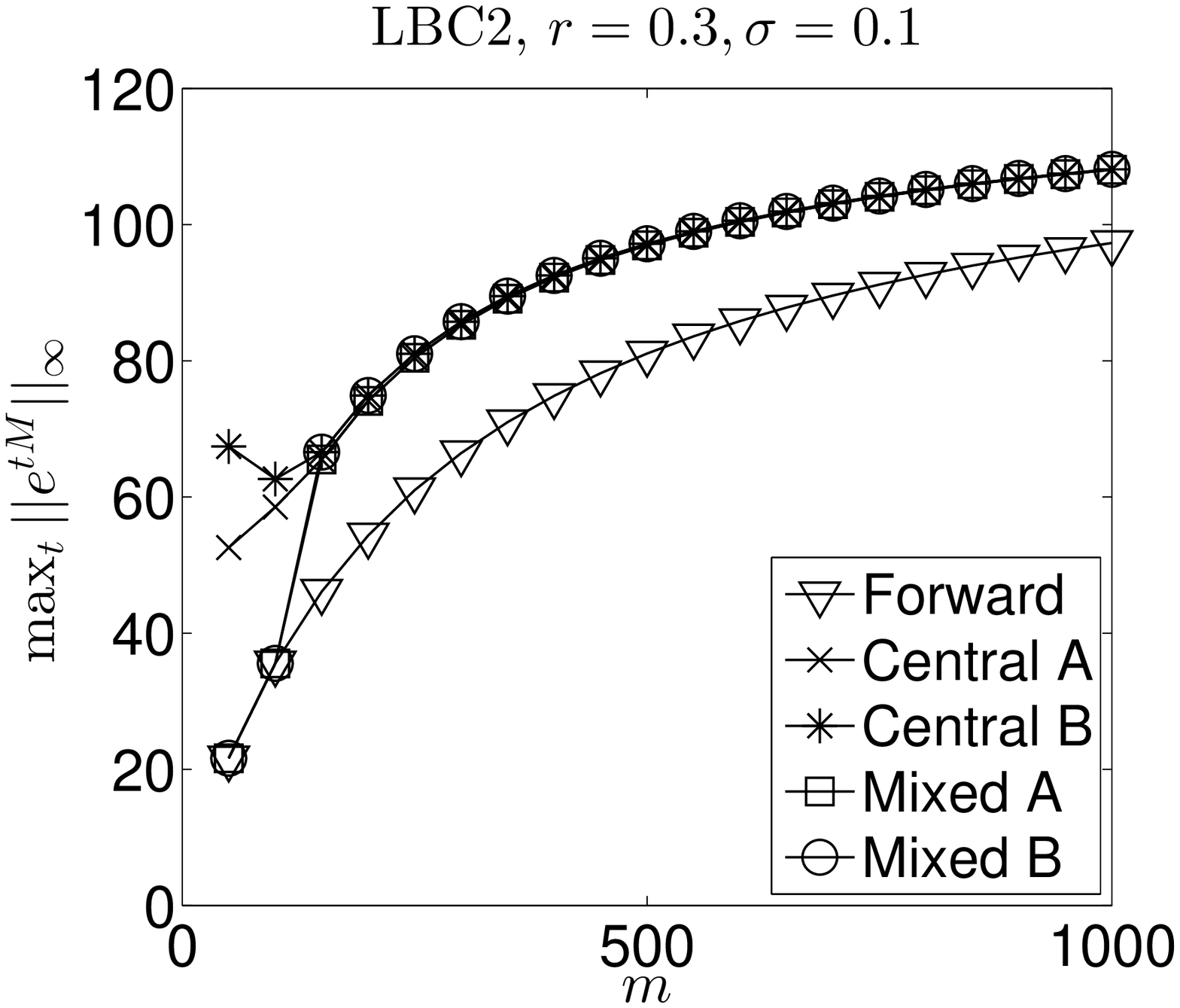}\\
\includegraphics[width=0.5\textwidth]{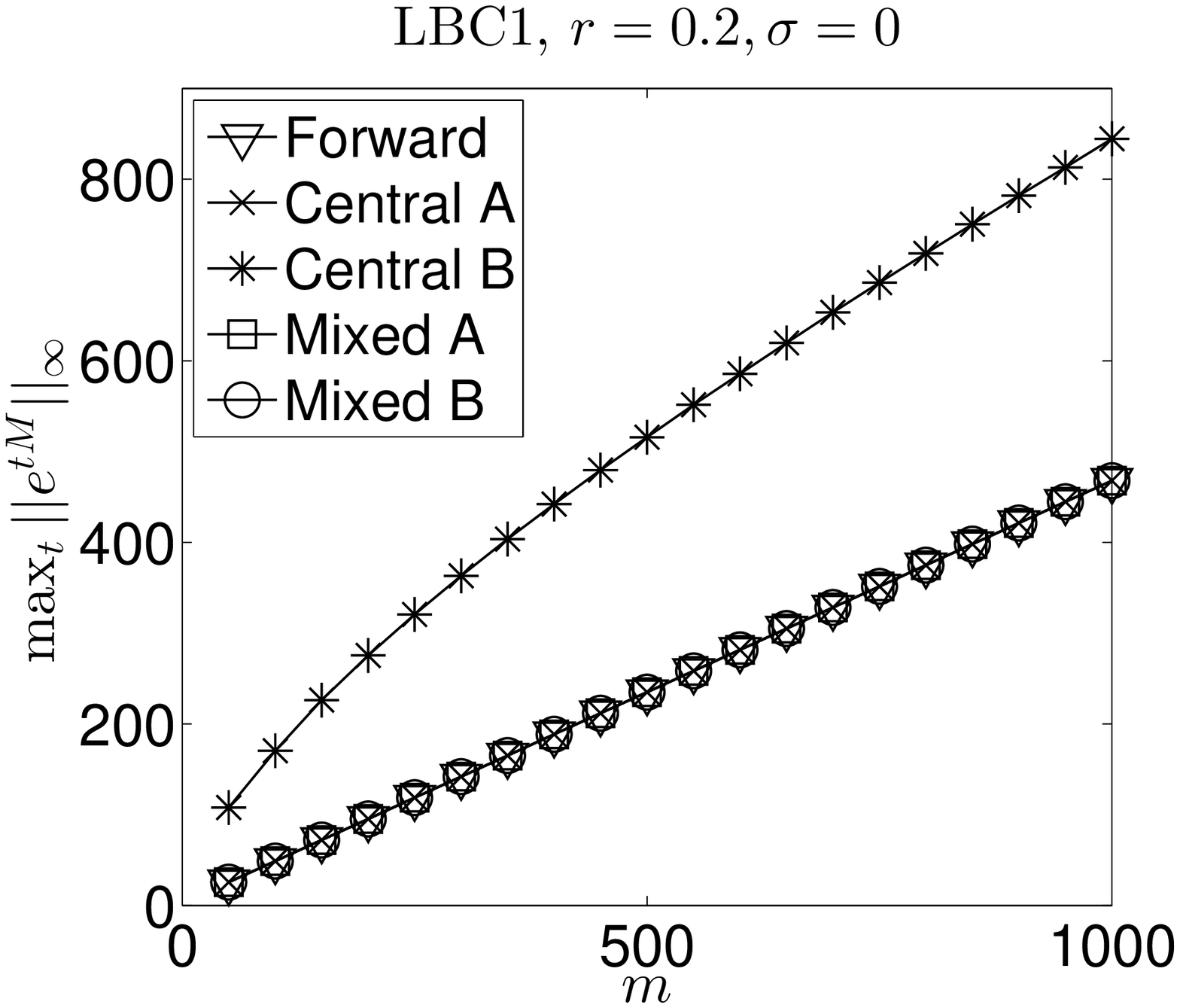}
&\includegraphics[width=0.5\textwidth]{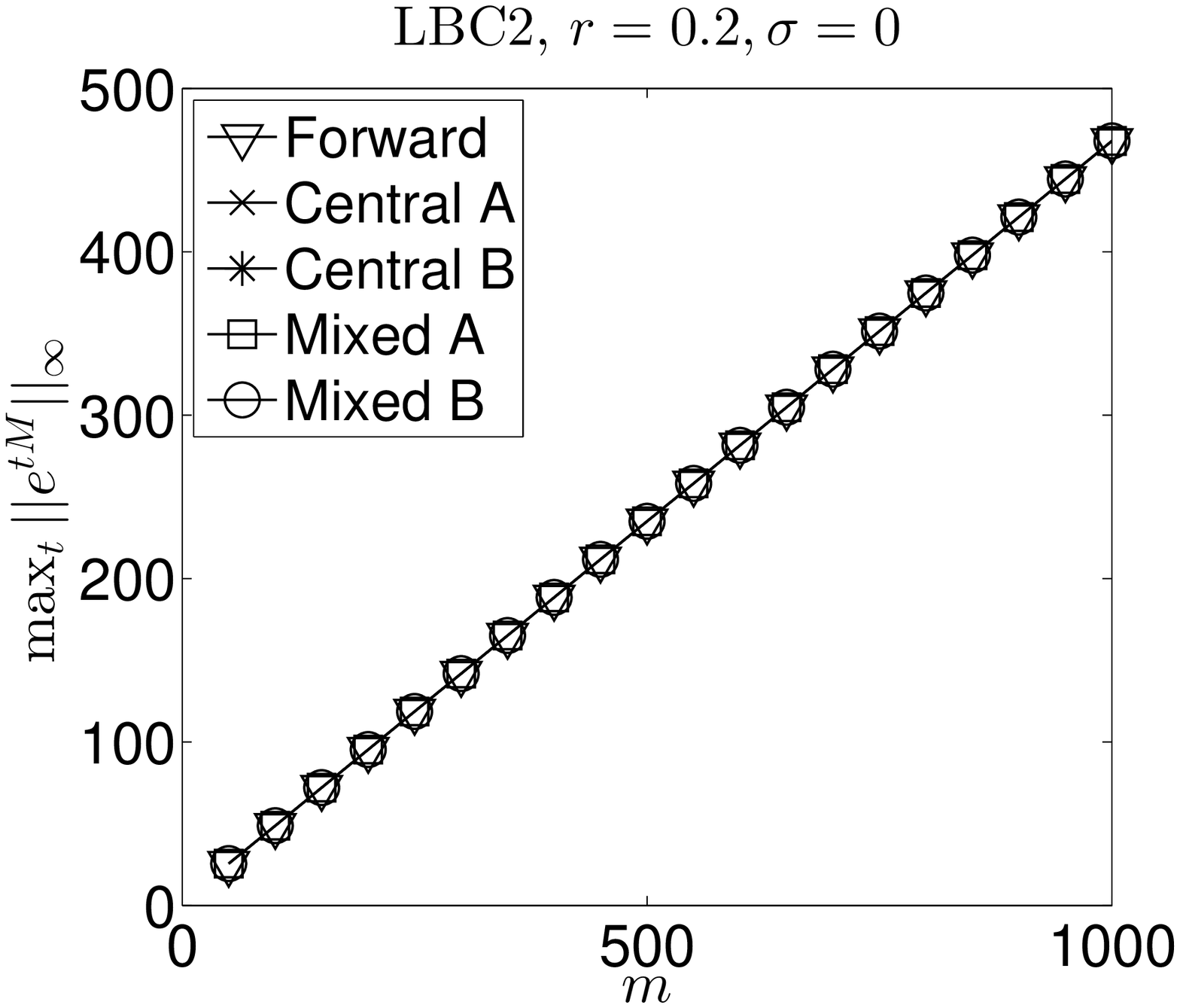}
\end{tabular}
\end{center}
\caption{Plot of $\max\{||e^{tM}||_\infty:t=0,1,\ldots,100\}$ vs.
$m=50,100,\ldots,1000$.
On the top row~$r=0.1$ and $\sigma=0.3$; on the middle row~$r=0.3$
and $\sigma=0.1$; on the bottom row~$r=0.2$ and $\sigma=0$.
The left column concerns the LBC1 discretization of the
linear boundary condition; the right column concerns
the LBC2 discretization.
In all cases $E=100$ and $S=400$.
Notice the different scales on the vertical axes.}
\label{fig:Stability}
\end{figure}

% new
First,
the stability of the ten FD discretizations discussed
above is numerically investigated.
For each FD discretization we consider the maximum of
$||e^{tM}||_\infty$ over $t\geq 0$ for the dimensions
$m=50,100,\ldots,1000$.
We employed the Matlab function {\tt expm} to evaluate the
matrix exponential and computed the values of $||e^{tM}||_\infty$
for $t=0,1,\ldots,100$; this was found to give a reliable
estimate\footnote{Except for Central A and B with LBC2 if
$\sigma=0$, see the discussion on this case in the text
further on.} for the maximum over all $t\ge 0$.
For the experiments three pairs $r,\sigma$ have been chosen.
Figure~\ref{fig:Stability} displays the obtained results.
On the top row~$r=0.1$ and $\sigma=0.3$; on the middle
row~$r=0.3$ and $\sigma=0.1$; on the bottom row~$r=0.2$
and $\sigma=0$.
The left column represents the LBC1 treatment of the linear
boundary condition and the right column represents the LBC2
treatment.
% new
Note
that the scales on the vertical axes vary.

% note: maximum usually attained at t=T.

If $r=0.1$ and $\sigma=0.3$, then the Mixed A and B
discretizations are identical to Central A and B,
respectively; in this case the conditions (\ref{sufA}),
(\ref{sufB}) always hold.
On the other hand, if $r=0.2$ and $\sigma=0$, then
Mixed A and B both reduce to the Forward discretization;
in this case (\ref{sufA}), (\ref{sufB}) never hold.
Finally, if $r=0.3$ and $\sigma=0.1$, then Mixed A
and B form an actual mix of Forward with Central A
and B, respectively.

%\begin{table}
%   \centering
%       \begin{tabular}{c|c|c}
%           $m$&Fraction Mixed A&Fraction Mixed B\\
%           \hline
%           50&60.0\%&64.0\%\\
%           100&28.0\%&30.0\%\\
%           150&13.3\%&13.3\%\\
%           200&11.0\%&10.5\%\\
%           250&9.2\%&8.8\%\\
%           300&8.0\%&7.7\%\\
%           350&6.9\%&6.9\%\\
%           400&6.3\%&6.0\%\\
%           450&5.6\%&5.6\%\\
%           500&5.2\%&5.0\%\\
%           550&4.7\%&4.7\%\\
%           600&4.3\%&4.3\%\\
%           650&4.0\%&4.0\%\\
%           700&3.9\%&3.7\%\\
%           750&3.6\%&3.6\%\\
%           800&3.4\%&3.4\%\\
%           850&3.2\%&3.2\%\\
%           900&3.0\%&3.0\%\\
%           950&2.8\%&2.8\%\\
%           1000&2.7\%&2.7\%
%       \end{tabular}
%   \caption{The fraction of grid points where the first order forward scheme
%   is used in Mixed A (column 2) and in Mixed B (column 3) vs. the dimension $m$ (column 1)
%   for the case \mbox{$u_{ss}(s_{m+1})=u_{ss}(S,t)=0$} and $r=0.3,\sigma=0.1$.}
%   \label{tab: Fraction mixed A en B stab}
%\end{table}

Considering the left column of Figure~\ref{fig:Stability},
a main observation is that the results for the Forward and the
Mixed A and B discretizations with the LBC1 treatment agree
with the stability bound (\ref{bound}) given by
Theorem~\ref{general thm}.
Indeed, for all  $r$, $\sigma$ pairs the maximum of
$||e^{tM}||_\infty$ over $t$ is
% new
found to be
directly proportional to $m$, and moreover,
\begin{equation}\label{expMnumer}
\max_{t\geq 0} ||e^{tM}||_\infty\approx \frac{2S}{h_{m+2}}\,.
\end{equation}
If $r=0.1$ and $\sigma=0.3$, then the results for Central~A
and B with LBC1 clearly also agree with (\ref{bound}), as
% new
these discretizations
coincide with Mixed~A and B.
We do not have a theoretical stability result for the central
discretizations for the other two $r$,~$\sigma$ pairs, but
it is interesting that the estimate (\ref{expMnumer}) is also
observed in the left column for Central~A and B if $r=0.3$
and $\sigma=0.1$, and for Central~A if $r=0.2$ and $\sigma=0$.
Although (\ref{expMnumer}) is not found for Central~B in
the latter case, the maximum of $||e^{tM}||_\infty$ over
$t$ still appears to be at most directly proportional
to~$m$.

Considering the right column of Figure~\ref{fig:Stability},
we observe in the top and middle rows that for all five FD
discretizations with the LBC2 treatment there appears to
be an upper bound on $||e^{tM}||_\infty$ that is uniform
both in $t$ and the dimension~$m$.
This upper bound is small in the top row and larger in
the middle row.
It appears to increase if the ratio $\sigma^2/r$ decreases.
In the limit case of $\sigma = 0$, displayed in the bottom
row, one obtains again linear growth with~$m$ for the
Forward and Mixed A and B discretizations with LBC2.
In fact, it is readily seen that
% new
these
three discretizations
% new
then all reduce to Forward with LBC1.
We notice that similar observations as here concerning
LBC2 were made in \cite{WFV04}.

If $\sigma = 0$, then the results for the Central A and B
discretizations with LBC2 are not displayed in the figure.
For each considered dimension $m$, the obtained maximum
of $||e^{tM}||_\infty$ is very large in this case.
For example, if $m=100$, then the maximum values (over
$t=0,1,\ldots,100$) for Central A and B are equal to
\mbox{$2.5\times 10^3$ and $1.5\times 10^5$}, respectively,
and if $m=200$, then they are equal to
\mbox{$9.8\times 10^3$ and $5.8\times 10^5$}.
In fact, there appears to be growth with $m^2$.
%We mention that theoretical stability results for the
%LBC2 treatment are

When $\sigma$ is not small, then comparing the results in
the left and right columns of Figure~\ref{fig:Stability}
one may be inclined to prefer the LBC2 treatment over
LBC1: for both $r$, $\sigma$ pairs with nonzero $\sigma$
one observes, for each given FD discretization, that the
maximum of $||e^{tM}||_\infty$ is more favorable for LBC2
compared to LBC1.
However, considering the actual convergence behavior of
the FD discretizations, we find no essential difference
between the two treatments, as is illustrated next.

\subsection{Convergence experiments}
% new
Here
we numerically examine the convergence behavior
of the ten FD discretizations of the Black--Scholes PDE with
linear boundary condition discussed above.
For the experiments we consider a European call option with
exercise price $E$, so that
\begin{equation*}\label{IC}
u(s,0) = \max(0,s-E) \quad (s\ge 0),~~~
u(0,t) = 0 \quad (0\le t \le T),
\end{equation*}
and set $E=100$, $T=5$, $S=2000$.
In the experiments we compute the~maximum norm of the
spatial discretization errors
\[
e_h(T) = u_h(T)-U(T)
\]
for a sequence of values $m$ with $10^2 \le m \le 10^4$.
Here $u_h$ denotes the restriction to the spatial grid
of the exact call option price function, given by the
Black--Scholes formula.
The semidiscrete solution vector $U(T)$ to (\ref{ODE}),
(\ref{general M}) is approximated, with sufficiently
high accuracy, by applying the Crank--Nicolson method
using $N=10^4$ time steps.
In view of the nonsmooth initial condition, two initial
(damping) substeps are taken with the implicit Euler
method.
This approach is often referred to as Rannacher time
stepping.

%We determine the exact solution of the Black--Scholes
%equation with the Matlab-function \texttt{blsprice}.

\begin{figure}
\begin{center}
\begin{tabular}{cc}
\includegraphics[width=0.5\textwidth]{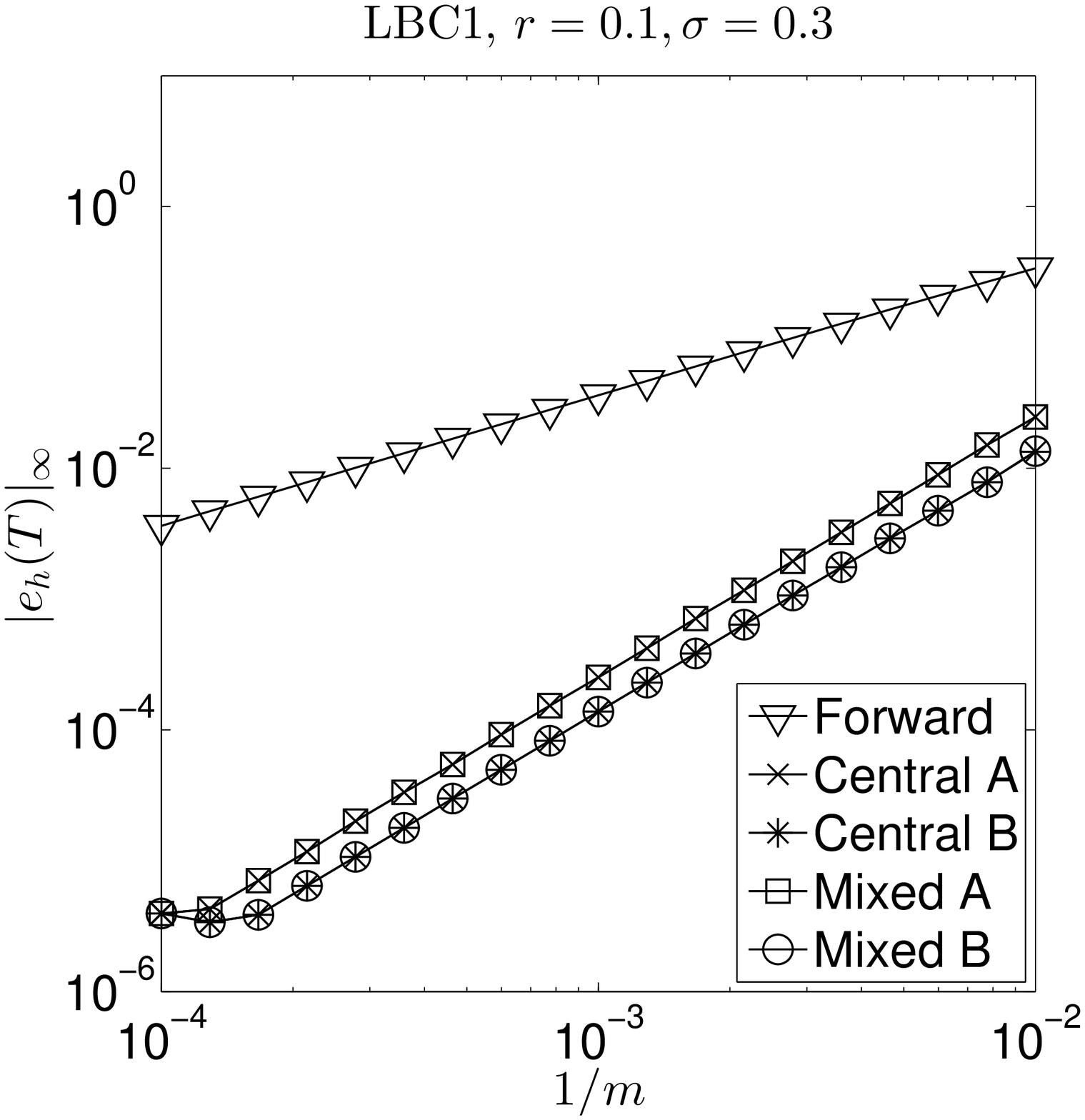}
&\includegraphics[width=0.5\textwidth]{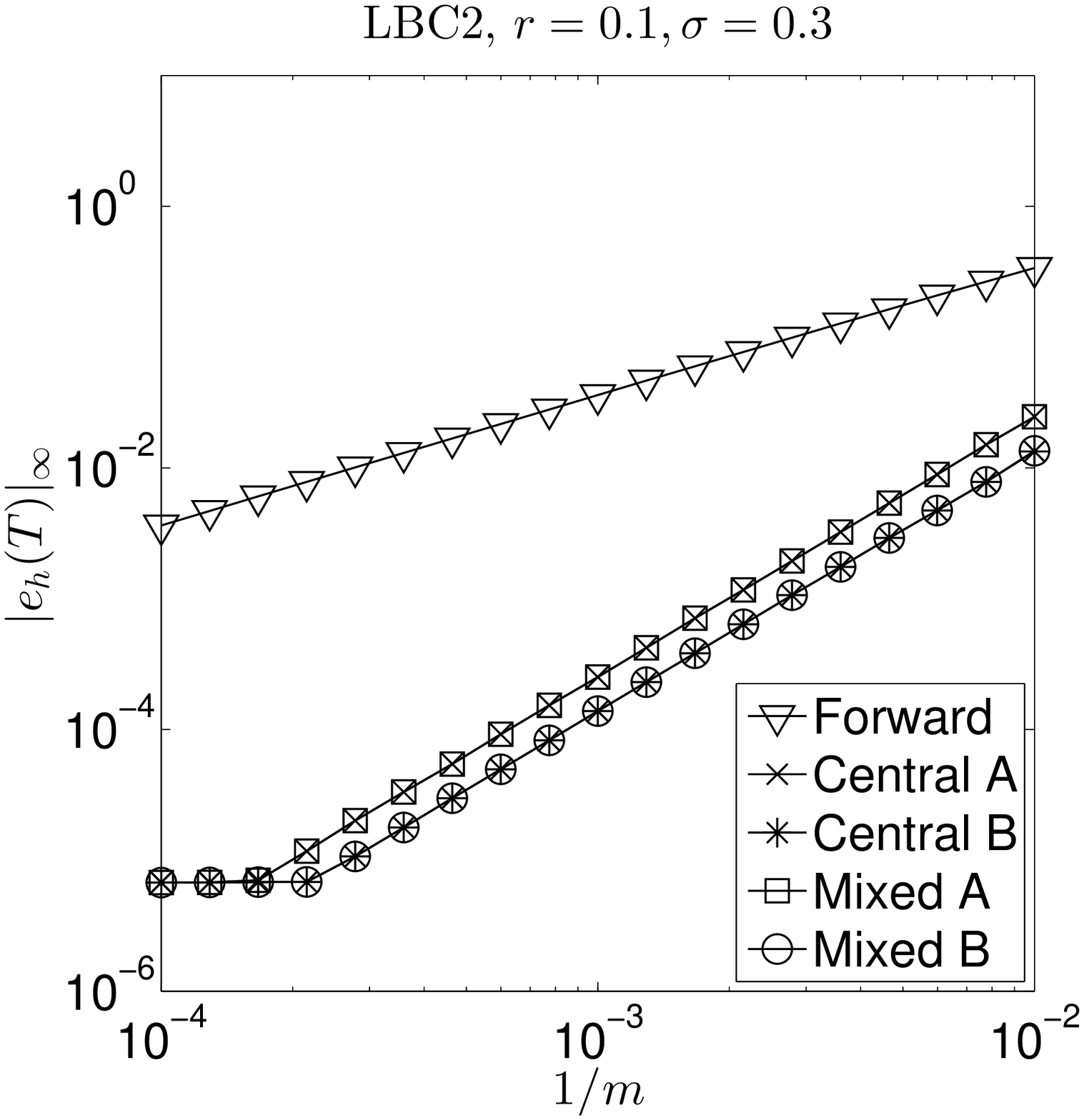}\\
\includegraphics[width=0.5\textwidth]{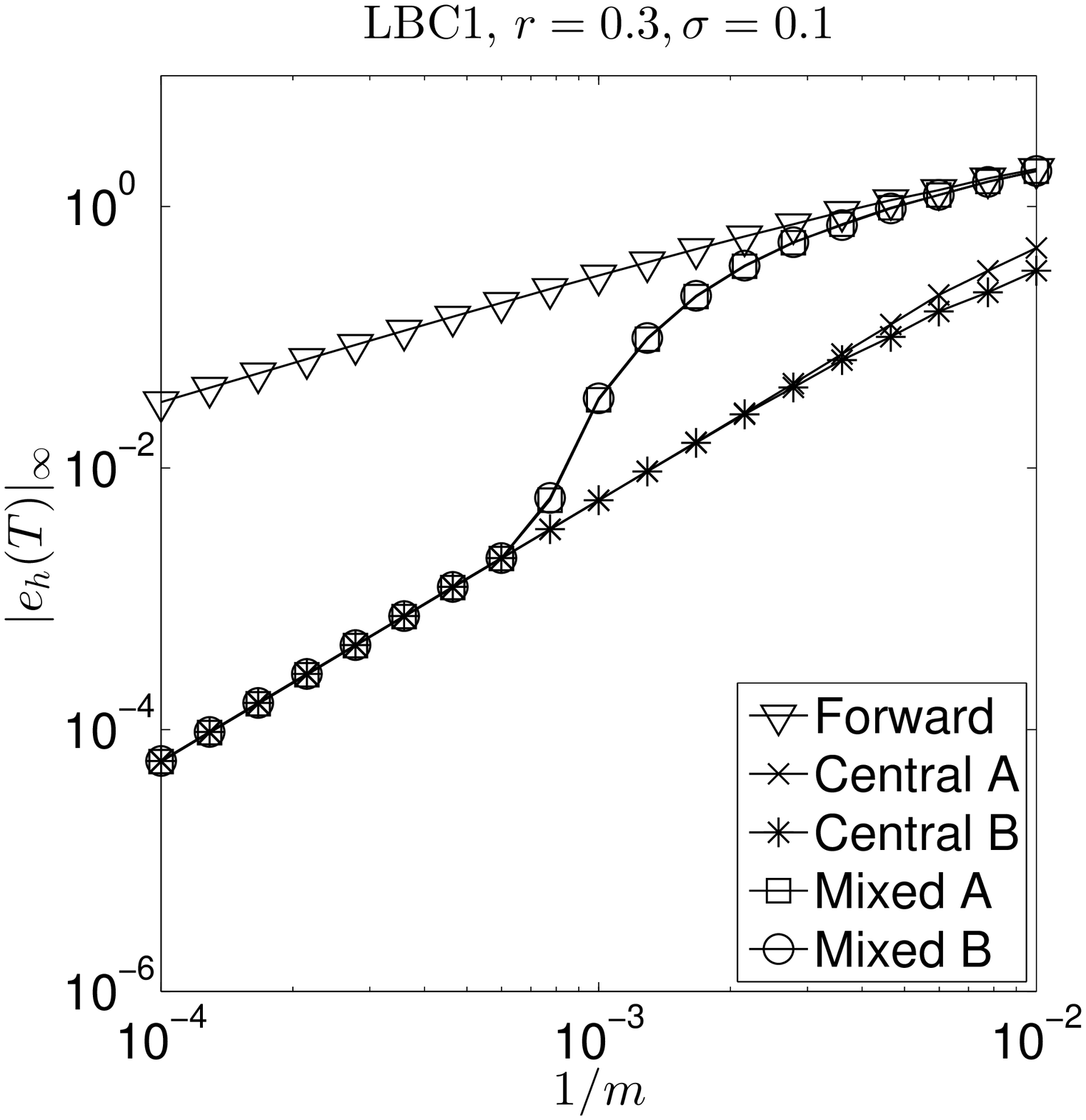}
&\includegraphics[width=0.5\textwidth]{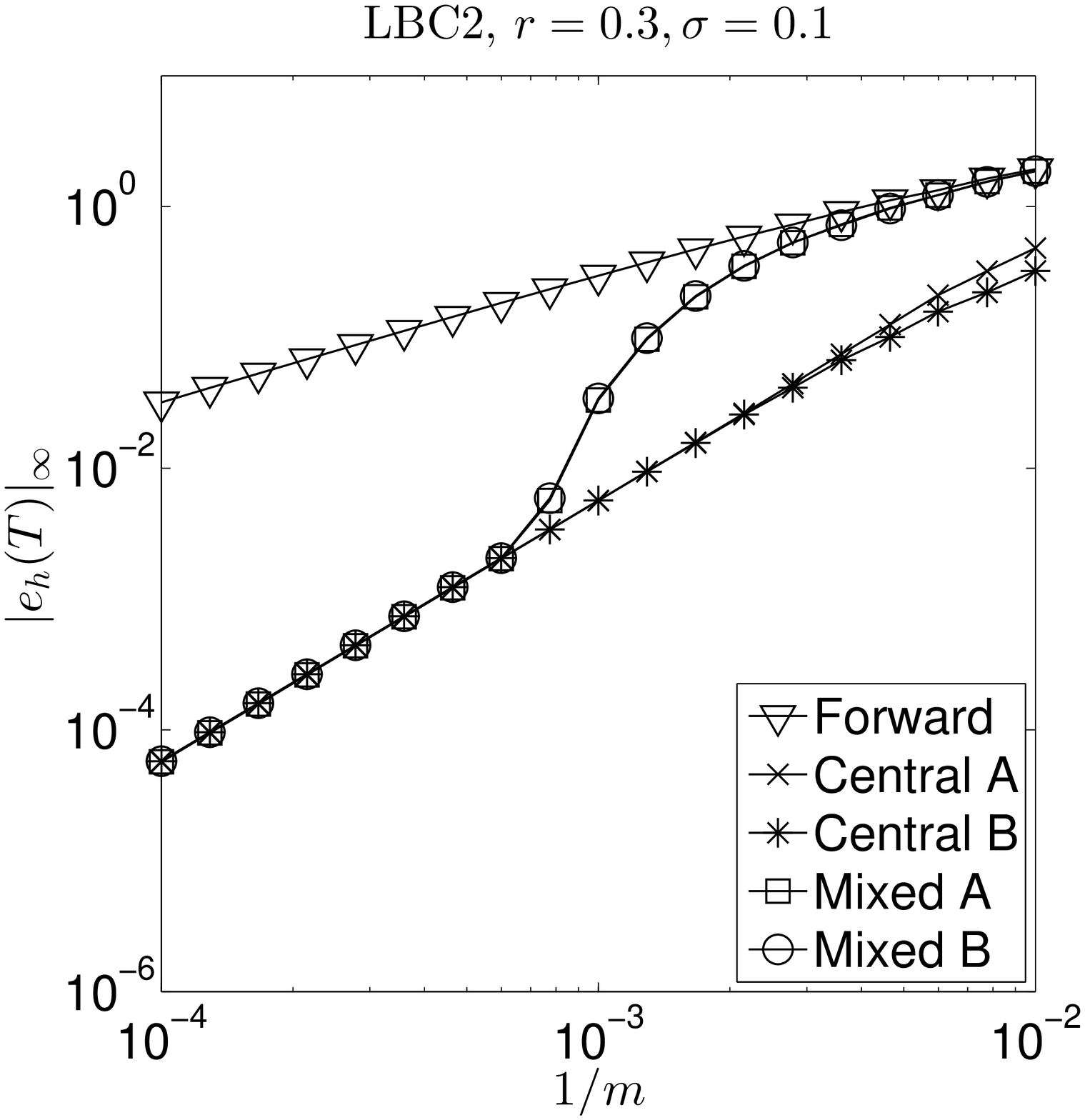}
\end{tabular}
\end{center}
\caption{Plot of $|e_h(T)|_\infty$ vs. $1/m$\, for values
$10^2\leq m\leq 10^4$ equally spaced on a logarithmic scale.
On the top row~$r=0.1$ and $\sigma=0.3$; on the bottom row
$r=0.3$ and $\sigma=0.1$.
The left column concerns the LBC1 discretization of the
linear boundary condition; the right column concerns
the LBC2 discretization.
In all cases $E=100$, $T=5$, $S=2000$.
}\label{fig:Convergence}
\end{figure}

The results are displayed in Figure \ref{fig:Convergence}.
On the first row $r=0.1$, $\sigma=0.3$ and on the second
row $r=0.3$, $\sigma=0.1$.
The left column represents the LBC1 treatment of the linear
boundary condition and the right column represents the LBC2
treatment.

As a first observation it is clear that LBC1 and LBC2
almost always lead to the same spatial discretization
error in the experiments.

If $r=0.1$ and $\sigma=0.3$, then the Mixed A and B
discretizations are identical to Central A and B,
respectively.
For these four discretizations we find, by least-squares
approximation in the region $m\le 5000$, an order of
convergence equal to 2.0.
Once $m$ gets larger than 5000, then the (fixed) error
due to the linear boundary condition dominates.
The Forward discretization has an observed order of
convergence equal to 1.0.

If $r=0.3$ and $\sigma=0.1$
% new
the observed orders of convergence are 0.9 for Forward,
2.0 for Central A and 1.9 for Central B.
The Mixed A and B discretizations are
% new
in this case an
actual mix of Forward with Central A and B, respectively.
Table~\ref{tab: Fraction mixed A en B} gives for each~$m$
under consideration the fraction of grid points where the
Forward discretization is used.
Clearly, for small $m$ this fraction is large and Mixed A
and B
% new
show spatial discretization errors similar to
% new
those for Forward, whereas for
large $m$ this fraction is small and Mixed A and B
% new
show
errors similar to
% new
those
for the Central A and B discretizations.
A rigorous analysis of this correlation
% new
will be left for future research.
\\

% new
\begin{table}[h!]
    \centering
        \begin{tabular}{c|c|c}
            $m$&Fraction Mixed A&Fraction Mixed B\\
            \hline
            100&57.0\%&58.0\%\\
            129&47.3\%&48.8\%\\
            167&34.1\%&36.5\%\\
            215&9.8\%&11.2\%\\
            278&7.9\%&7.6\%\\
            359&6.4\%&6.4\%\\
            464&5.2\%&5.2\%\\
            599&4.2\%&4.2\%\\
            774&3.4\%&3.4\%\\
            1000&2.7\%&2.7\%\\
            1292&2.1\%&2.1\%\\
            1668&1.7\%&1.7\%\\
            2154&1.3\%&1.3\%\\
            2783&1.0\%&1.0\%\\
            3594&0.8\%&0.8\%\\
            4642&0.6\%&0.6\%\\
            5995&0.5\%&0.5\%\\
            7743&0.4\%&0.4\%\\
            10000&0.3\%&0.3\%
        \end{tabular}
\caption{Fraction of grid points where the Forward discretization
is used in Mixed A and Mixed B vs. $m$\, in the case of $r=0.3$,
$\sigma=0.1$ and LBC1.}
\label{tab: Fraction mixed A en B}
\end{table}

\vfill\eject
%%%%%%%%%%%%%%%%%%%%%%%%%%%%%%%%%%%%%%%%%%%%%%%%%%%%%%%%%%%%%%%%%%%%%%%%%%%%%%%%%%%%%%
\section{Time discretization}\label{timediscr}
\setcounter{equation}{0}
\setcounter{theorem}{0}

In this section we study the time discretization of
the semidiscrete system (\ref{ODE}), (\ref{general M})
by the well-known family of $\theta$-methods, which
includes the popular Crank--Nicolson method (trapezoidal
rule) and implicit Euler method as special cases.
As noted in Section~\ref{intro}, the subsystem of
ODEs involving the matrix $C$ could be solved exactly,
but it is more interesting and useful, both from a
theoretical and practical point of view, to consider
the time discretization of the semidiscrete system
(\ref{ODE}), (\ref{general M}) as a
whole.
% new
%\footnote{It may appear that the stability and
%convergence analysis simplifies much if the subsystem
%involving $C$ is solved exactly, but it turns out that
%this is not the case.}

Let parameter $\theta \in [\tfrac{1}{2},1]$ be given
and fixed.
Let step size $\Delta t = T/N$ with integer $N\ge 1$ be
given and define $t_n = n\, \Delta t$, $b_n = b(t_n)$
for $n=0,1,2,\ldots,N$.
The {\it $\theta$-method} generates, in a successive way,
for $n=1,2,\ldots,N$ an approximation $U_n$ to $U(t_n)$ by
\[
U_n = U_{n-1} + (1-\theta)\,\Delta t\, (MU_{n-1} + b_{n-1})
+ \theta\, \Delta t\, (MU_{n} + b_{n}).
\]
The choices $\theta = \frac{1}{2}$ and $\theta = 1$ yield,
respectively, the Crank--Nicolson method and implicit
Euler method.
The above recurrence relation can be written as
\begin{equation}\label{theta}
U_n = \varphi (\Delta t\, M) U_{n-1} + (I-\theta \Delta t\, M)^{-1}
[(1-\theta)\,\Delta t\, b_{n-1} + \theta\, \Delta t\, b_{n}],
\end{equation}
where $\varphi$ is the so-called stability function of the method,
given by
\[
\varphi (z) = \frac{1+(1-\theta)z}{1-\theta z} \quad ({\rm for}~ z\in\C),
\]
and
\[
\varphi (X) =
(I-\theta X)^{-1} (I+(1-\theta)X) =
(I+(1-\theta)X)(I-\theta X)^{-1}
\]
\vskip0.2cm\noindent
for square matrices $X$ such that $I-\theta X$ is invertible.
% with $I$ the identity matrix of the same size as $X$.

We first study the stability of the fully discrete process
(\ref{theta}).
The
% new
subsequent two lemmas can be viewed as analogues of Lemmas
\ref{lemma expM}, \ref{lemma expC} for the semidiscrete
system.

\begin{lemma}\label{lemma varphiM}
For $1\le n\le N$ there holds
\[
\varphi (\Delta t\, M)^n =
\left(\begin{array}{c|c}
\varphi (\Delta t\, A)^n&Y_n\\
\hline
O&\varphi (\Delta t\, C)^n
\end{array}\right)
\]
where
\begin{subeqnarray}\label{Ymatrix}
Y_1 &=& (I-\theta \Delta t\, A)^{-1} \Delta t\, B (I-\theta \Delta t\, C)^{-1},\\
Y_n &=& \sum_{j=0}^{n-1}\, \varphi (\Delta t\, A)^{n-j-1}\, Y_1\, \varphi (\Delta t\, C)^{j}.
\end{subeqnarray}
\end{lemma}
\textbf{Proof}~
The formula is readily obtained by induction to $n$ and noting that
\begin{equation}\label{invthetaM}
(I-\theta \Delta t\, M)^{-1} =
\left(\begin{array}{c|c}
(I-\theta \Delta t\, A)^{-1}&\theta\, Y_1\\
\hline
O&(I-\theta \Delta t\, C)^{-1}
\end{array}\right).
\end{equation}
\begin{flushright}
$\Box$
\end{flushright}\vskip0.2cm

\begin{lemma}\label{lemma varphiC}
For $1\le n\le N$ there holds
\begin{equation*}\label{varphiC}
||\varphi (\Delta t\, C)^n||_\infty =
x^n+\left(1-x^n \right)\frac{2S}{h_{m+2}} \quad
{\it with}~~ x = \varphi(-r \Delta t).
\end{equation*}
\end{lemma}
\textbf{Proof}~With the eigendecomposition of $C$ it is easily
verified that
\begin{equation*}\label{varphiCform}
\varphi (\Delta t\, C)^n = \frac{1}{h_{m+2}}\left(\begin{array}{cc}
Sx^n-s_{m+1}&s_{m+1}(1-x^n)\\
S(x^n-1)&S-s_{m+1}x^n
\end{array}\right)
\end{equation*}
and the rest of the proof is similar to that of
Lemma \ref{lemma expC}, using
% new
%that
$|x|\le 1$.
\begin{flushright}
$\Box$
\end{flushright}\vskip0.2cm

Concerning the discretization on the spatial domain
$[s_1, s_m]$ we shall
% new
assume
in the following that there exists a real constant $K$,
independent of the dimension $m$ and number of
time steps $n$, such that
\begin{equation}\label{stabconst}
||\varphi (\Delta t\, A)^n||_\infty \le K \quad
{\rm whenever}~~\Delta t =T/N,~ 0\le n \le N,~ N\ge 1.
\end{equation}
In the literature much attention has been paid to
establishing (\ref{stabconst}), under a variety
of conditions on the matrix $A$.
For the implicit Euler method ($\theta =1$) the neat
result is well-known that (\ref{stabconst}) is
fulfilled with $K=1$ whenever $\mu_\infty [A] \le 0$,
% new
cf.~e.g.~\cite{HW02,S83}.
% new
Hence, this is guaranteed under the condition (\ref{condition}).
For all other time discretization methods, however,
the available results in the literature
% new
%that
implying
(\ref{stabconst}) with a constant $K$
independent of $m$ and $n$ require, to the best
of our knowledge, stronger conditions on $A$.
Notably, conditions on the resolvent, the numerical
range and
% new
%the
pseudospectra have been extensively
investigated in the literature,
% new
cf.~e.g.~\cite{S98,TE05}.
Some of these results appear to be useful in our
current application,
% new
%of the semidiscretized Black--Scholes PDE,
but a verification of the
pertinent
% new
conditions
on $A$
% new
is highly non-trivial.
% new
As our main interest in this paper lies in
studying (the implications of) the discretized
linear boundary condition
% new
on $[s_{m+1},s_{m+2}]$,
which corresponds to the matrices
$B$ and $C$, we shall leave the analysis of
(\ref{stabconst}) when
% new
$\tfrac{1}{2} \le \theta < 1$ for future
research.\\

% new
The next theorem can be regarded as a discrete analogue
to Theorem \ref{general thm}.

\begin{theorem}\label{general thm_varphi}
If\, {\rm (\ref{condition})} and {\rm (\ref{stabconst})}
then for $1\le n\le N$,
\begin{equation*}\label{bound varphiM}
x^n+\left(1-x^n \right)\frac{2S}{h_{m+2}}
\, \leq \, ||\varphi (\Delta t\, M)^n||_\infty \, \leq \,
K+\frac{4(K+1)S}{h_{m+2}}
\end{equation*}
with $x = \varphi(-r \Delta t)$.
\end{theorem}
\textbf{Proof}~
% new
For any integer $j\ge 0$, let the rational
function $\psi_j$ be defined by
\begin{equation*}\label{psi_j}
\psi_j(z) = \frac{\varphi(z)^j}{1-\theta z}\quad (z\in\C).
\end{equation*}
Consider the formula for $\varphi (\Delta t\, M)^n$
given by Lemma \ref{lemma varphiM}.
The lower bound on its maximum norm is clear by
Lemma \ref{lemma varphiC}.
To prove the upper bound, write
\begin{eqnarray*}
  Y_n &=& \sum_{j=0}^{n-1}\,
\varphi (\Delta t\, A)^{n-j-1}\, (I-\theta \Delta t\, A)^{-1}
\Delta t\, B\, (I-\theta \Delta t\, C)^{-1} \varphi (\Delta t\, C)^{j} \\
   &=& \sum_{j=0}^{n-1}\,
\varphi (\Delta t\, A)^{n-j-1}\, (I-\theta \Delta t\, A)^{-1}
\Delta t\, B\, \psi_j (\Delta t\, C).
\end{eqnarray*}
It holds that
\begin{equation}\label{BpsiC}
B\, \psi_j (\Delta t\, C) =
\frac{\gamma_{m}}{h_{m+2}}
\left[\begin{array}{cc}(Sx_j-s_{m+1})e_m&
(s_{m+1}-s_{m+1}x_j)e_m\end{array}\right]
\end{equation}
where $x_j = \psi_j(-r \Delta t)\in [-1,1]$.
Both columns of this matrix are of the form
\[
f_j = (\phi_0 + \phi_1 x_j) e_m
\]
with real numbers $\phi_0$, $\phi_1$ independent of $j$.

We have
\begin{eqnarray*}
&&   \Delta t\, \sum_{j=0}^{n-1}\, \varphi (\Delta t\, A)^{n-j-1}\, (I-\theta \Delta t\, A)^{-1} \\
&& = \Delta t\, \sum_{k=0}^{n-1}\, \varphi (\Delta t\, A)^{k}\, (I-\theta \Delta t\, A)^{-1} \\
&& = \Delta t\, (\varphi (\Delta t\, A)^n - I)(\varphi (\Delta t\, A) - I)^{-1} (I-\theta \Delta t\, A)^{-1} \\\\
&& = \Delta t\, (\varphi (\Delta t\, A)^n - I)(I+(1-\theta)\Delta t\, A - I + \theta \Delta t\, A)^{-1} \\\\
&& = (\varphi (\Delta t\, A)^n - I) A^{-1}.
\end{eqnarray*}
By similar algebraic manipulations, there follows
\begin{eqnarray*}
&&   \Delta t\, \sum_{j=0}^{n-1}\, \varphi (\Delta t\, A)^{n-j-1}\, (I-\theta \Delta t\, A)^{-1} x_j \\
&& = ( \varphi (\Delta t\, A)^n - \varphi(-r \Delta t)^n I ) (rI+A)^{-1}.
\end{eqnarray*}
Consequently,
\begin{eqnarray}\label{sumformula}
&&   \Delta t\, \sum_{j=0}^{n-1}\, \varphi (\Delta t\, A)^{n-j-1}\, (I-\theta \Delta t\, A)^{-1} f_j \\
&& = \phi_0\, (\varphi (\Delta t\, A)^n - I) A^{-1}e_m +
     \phi_1\, ( \varphi (\Delta t\, A)^n - \varphi(-r \Delta t)^n I ) (rI+A)^{-1}e_m.\nonumber
\end{eqnarray}
\vskip0.2cm\noindent
Continuing from here along the same lines as in the proof of
Theorem \ref{general thm}, and using the condition (\ref{stabconst}),
we arrive at
\[
||Y_n||_\infty \le \frac{4(K+1)S}{h_{m+2}}\,.
\]
Together with Lemma  \ref{lemma varphiC}, the stated upper
bound on $||\varphi (\Delta t\, M)^n||_\infty$ now directly
follows.
\begin{flushright}
$\Box$
\end{flushright}\vskip0.2cm

In the subsequent convergence analysis of the process
(\ref{theta}), the matrix
\[
\psi_n(\Delta t\, M) =
\varphi (\Delta t\, M)^n (I-\theta \Delta t\, M)^{-1}
\]
arises.
By Lemma \ref{lemma varphiM} and formula (\ref{invthetaM}),
\begin{equation*}\label{varpsiM}
\psi_n(\Delta t\, M) =
\left(\begin{array}{c|c}
\psi_n(\Delta t\, A)&
\theta\, \varphi (\Delta t\, A)^n\,  Y_1 +Y_n\,(I-\theta \Delta t\, C)^{-1}\\
\hline
O&\psi_n(\Delta t\, C)
\end{array}\right).
\end{equation*}
For the analysis below we need upper bounds on
the maximum norms of the constituent submatrices.
Putting $Y_0=O$, there holds

\begin{lemma}\label{lemma psiM}
If\, {\rm (\ref{condition})} and {\rm (\ref{stabconst})}
then for $0\le n\le N-1$,
\begin{subeqnarray}\label{psiM}
||\psi_n(\Delta t\, A)||_\infty &\le& K\,,\phantom{\frac{4S}{h_{m+2}}}\\
||\psi_n(\Delta t\, C)||_\infty &\le& \frac{4S}{h_{m+2}}\,,\\
||\varphi (\Delta t\, A)^n\, Y_1||_\infty &\le& \frac{4KS}{h_{m+2}}\,,\\
||Y_n\, (I-\theta \Delta t\, C)^{-1}||_\infty &\le& \frac{4(K+1)S}{h_{m+2}}\,.
\end{subeqnarray}
\end{lemma}
\textbf{Proof}~The bound (\ref{psiM}a) follows directly from
(\ref{stabconst}) and the fact that (\ref{condition}) implies
$||(I-\theta \Delta t\, A)^{-1}||_\infty \le 1$.
The bound (\ref{psiM}b) is obtained using the same arguments
as in the proof of Lemma \ref{lemma varphiC}.
In order to prove (\ref{psiM}c) we note that
\[
(I-\theta \Delta t\, A)^{-1} =
(\varphi (\Delta t\, A) - I)(\Delta t\, A)^{-1}
\]
and using this gives
\begin{equation}\label{varphiAnY}
\varphi (\Delta t\, A)^n\, Y_1 =
(\varphi (\Delta t\, A)^{n+1}-\varphi (\Delta t\, A)^n)\,
A^{-1} B\, (I-\theta \Delta t\, C)^{-1}.
\end{equation}

\noindent
By (\ref{BpsiC}) with $j=0$,
\[
B\, (I-\theta \Delta t\, C)^{-1} =
\frac{\gamma_{m}}{h_{m+2}}
\left[\begin{array}{cc}(Sx_0-s_{m+1})e_m&
(s_{m+1}-s_{m+1}x_0)e_m\end{array}\right]
\]
where $x_0 = (1+ \theta\, r \Delta t)^{-1}\in (0,1)$.
If $f_0$ represents any of the two columns of this matrix,
then by a same argument as in the proof of Theorem
\ref{general thm} there follows
\[
|A^{-1}f_0|_\infty \le \frac{S}{h_{m+2}}\,.
\]
Consequently,
\[
||(\varphi (\Delta t\, A)^{n+1}-\varphi (\Delta t\, A)^n)\,
A^{-1} B\, (I-\theta \Delta t\, C)^{-1}||_\infty
\le 2K\cdot 2 \cdot \frac{S}{h_{m+2}} = \frac{4KS}{h_{m+2}}\,. % new
\]
The proof of the bound (\ref{psiM}d) is identical to that for
$||Y_n||_\infty$ given above, except that $\psi_j(z)$ needs to
be replaced by $\psi_j(z)/(1-\theta z)$.
\begin{flushright}
$\Box$
\end{flushright}\vskip0.2cm

To prove the convergence result for the time discretization
process (\ref{theta}), we also need the following result.

\begin{lemma}\label{lemma psiMg}
Assume {\rm (\ref{condition})} and {\rm (\ref{stabconst})}
hold.
Let $g:[0,S] \rightarrow \R$ be any given~continuously
differentiable function and
\[
w = \left(
      \begin{array}{c}
        g(s_{m+1}) \\
        g(s_{m+2}) \\
      \end{array}
    \right) \in \R^2.
\]
Then there exists $\xi \in (s_{m+1},s_{m+2})$ such that
for all $0\le n\le N-1${\rm :}

\begin{subeqnarray}\label{psiMg}
|\psi_n(\Delta t\, C)\, w|_\infty &\le& |g(S)|+ 2S |g'(\xi)|,\\ \nonumber \\
|\varphi (\Delta t\, A)^n\, Y_1\, w|_\infty &\le& 2K (|g(S)|+ S |g'(\xi)|),\\ \nonumber \\
|Y_n\, (I-\theta \Delta t\, C)^{-1}\, w|_\infty &\le& (K+1)(|g(S)|+ 2S |g'(\xi)|).
\end{subeqnarray}
\end{lemma}
\vskip0.2cm\noindent
\textbf{Proof}~Write $s=s_{m+1}$ and $h=h_{m+2}$.
Let $\xi \in (s,S)$ be such
that
\[
g(s) = g(S) - h g'(\xi).
\]
Then the vector $w$ can be written as
\[
w = g(S)
    \left(
      \begin{array}{c}
        1 \\
        1 \\
      \end{array}
    \right)
   - h g'(\xi)
   \left(
      \begin{array}{c}
        1 \\
        0 \\
      \end{array}
    \right).
\]
For any rational function $\psi$ with $\psi(0)=1$ there
holds
\[
\psi (\Delta t\, C) =
\frac{1}{h}\left(\begin{array}{cc}
Sx-s&s(1-x)\\S(x-1)&S-sx
\end{array}\right)
\]
where $x=\psi(-r\Delta t)$.
Application of this matrix to $w$, in the above form,
readily yields
\[
\psi (\Delta t\, C) w =
\left(
  \begin{array}{c}
    g(S)x - g'(\xi)(Sx - s) \\
    g(S)x - g'(\xi)(Sx - S) \\
  \end{array}
\right).
\]
Observe the important fact that there is no factor $1/h$ present here.\\

\noindent
(a) The bound (\ref{psiMg}a) is obtained upon taking $\psi=\psi_n$
and using $|x|\le 1$.\\

\noindent
(b) By formula (\ref{varphiAnY}),
\[
|\varphi (\Delta t\, A)^n\, Y_1\,w|_\infty \le  % new
2K\cdot|A^{-1} B\, (I-\theta \Delta t\, C)^{-1}\,w|_\infty.
\]
Considering $\psi(z) = \psi_0(z) = (1-\theta z)^{-1}$ yields
\[
A^{-1} B\, (I-\theta \Delta t\, C)^{-1}\,w =
(g(S)x - g'(\xi)(Sx - s))\, \gamma_m\, A^{-1} e_m
\]
with $x=(1+ \theta\, r \Delta t)^{-1}$.
As in the proof of Theorem \ref{general thm} we have
\mbox{$|\gamma_m\, A^{-1} e_m|_\infty \le 1$} and,
together with $0<x<1$, there follows
\[
|A^{-1} B\, (I-\theta \Delta t\, C)^{-1}\,w|_\infty
\le |g(S)|+ S |g'(\xi)|,
\]
which completes the proof of (\ref{psiMg}b).\\

\noindent
(c) By formula (\ref{Ymatrix}),
\[
Y_n (I-\theta \Delta t\, C)^{-1}\, w = \sum_{j=0}^{n-1}\,
\varphi (\Delta t\, A)^{n-j-1}\, (I-\theta \Delta t\, A)^{-1}
\Delta t\, B\, {\widetilde\psi}_j (\Delta t\, C) w
\]
with ${\widetilde\psi}_j (z) =  \psi_j(z)/(1-\theta z)$.
Taking $\psi = {\widetilde\psi}_j$ in the above general formula,
we get
\[
B\, {\widetilde\psi}_j (\Delta t\, C) w = (\phi_0 + \phi_1 {\widetilde x}_j) e_m,
\]
where
\[
{\widetilde x}_j = {\widetilde\psi}_j (-r\Delta t)~~,~~
\phi_0 = \gamma_m\, g'(\xi)s~~,~~\phi_1 = \gamma_m\, (g(S)-g'(\xi)S).
\]
It is convenient to set $x_j = \psi_j(-r\Delta t)$ and
${\widetilde\phi}_1 = \phi_1/(1+ \theta\, r \Delta t)$.
Then $\phi_1 {\widetilde x}_j = {\widetilde\phi}_1 x_j$
and formula (\ref{sumformula}) directly gives
\begin{eqnarray*}
&&Y_n (I-\theta \Delta t\, C)^{-1}\, w \\
&& = \Delta t\, \sum_{j=0}^{n-1}\, \varphi (\Delta t\, A)^{n-j-1}\,
(I-\theta \Delta t\, A)^{-1} (\phi_0 + \phi_1 {\widetilde x}_j) e_m \\
&& = \phi_0\, (\varphi (\Delta t\, A)^n - I) A^{-1}e_m +
     {\widetilde\phi}_1\, ( \varphi (\Delta t\, A)^n - \varphi(-r \Delta t)^n I ) (rI+A)^{-1}e_m.
\end{eqnarray*}
% new
\newpage
Using that $|\gamma_m\, A^{-1} e_m|_\infty \le 1$ and
$|\gamma_m\, (rI+A)^{-1} e_m|_\infty \le 1$ yields
\begin{eqnarray*}
|Y_n\, (I-\theta \Delta t\, C)^{-1}\, w|_\infty
&\le&
(K+1)(|g'(\xi)|S + |g(S)- g'(\xi)S|) \\
&\le& (K+1)(|g(S)|+ 2S |g'(\xi)|),
\end{eqnarray*}
which proves the bound (\ref{psiMg}c).
\begin{flushright}
$\Box$
\end{flushright}
%\vskip0.2cm

% new
As in Section \ref{conv},
let the vector $u_h(t)$ be given by
\[
u_h(t)=(u(s_1,t),u(s_2,t),\ldots,u(s_{m+2},t))^{\rm T}\,,
\]
where $u$ is the exact solution to the initial-boundary
value problem for the Black--Scholes PDE (\ref{BS}) on
$0\le s\le S$ with linear boundary condition (\ref{LBC}).
The following theorem provides a useful estimate for
the \emph{space-time discretization error}, defined by
%for the $\theta$-method (\ref{theta}) by
\[
\eh_n = u_h(t_n) - U_n \quad (0\le n \le N).
\]
It essentially states that the estimate for the spatial
discretization error from Theorem \ref{convergence} remains
valid after time discretization up to a\, $c\cdot (\Delta t)^p$
term, where $p$ denotes the classical order of consistency
of the $\theta$-method and $c$ is a constant independent
of the spatial grid and the time step.

\begin{theorem}\label{convergence_time}
Let $p=1$ if\, $\tfrac{1}{2}< \theta \le 1$ and $p=2$ if\,
$\theta = \tfrac{1}{2}$.
Assume~that all partial derivatives of $u$ of orders $\le p+2$
exist and are continuous on $[0,S]\times [0,T]$.
Let $h^\ast> 0$ be given and let $\kappa$, $\eta$ be defined
by {\rm (\ref{kappa_eta})}.
Assume {\rm(\ref{condition})} and {\rm (\ref{stabconst})}~hold.
Then there exists a real constant $c$ (depending only on $u$,
$S$, $\theta$ and $K$) such that
\[
|\eh_n|_\infty \le  t_n \cdot
\left\{
K \cdot\max_{0\leq\vartheta\leq t_n} |\delta_h^L(\vartheta)|_\infty
+(K+\tfrac{1}{2})\kappa \cdot\max_{0\leq\vartheta\leq t_n}\eta(\vartheta)+
c\cdot (\Delta t)^{p}
\right\}
\]
whenever~ $0<h_{m+2} \le h^\ast$, $\Delta t =T/N$, $1\le n \le N$,
$N\ge 1$.
\end{theorem}
\textbf{Proof}~Let the local space-time error $\dh_n$ in the $n$-th
step of (\ref{theta}) be defined by
\begin{equation*}\label{deltahat}
u_h(t_n) = \varphi (\Delta t\, M) u_h(t_{n-1}) + (I-\theta \Delta t\, M)^{-1}
[(1-\theta)\,\Delta t\, b_{n-1} + \theta\, \Delta t\, b_{n}] + \dh_n.
\end{equation*}
Subtracting (\ref{theta}) from this yields
%(using $\eh_0 = 0$)
\begin{equation}\label{epshat}
\eh_n = \varphi (\Delta t\, M) \eh_{n-1} + \dh_n = \ldots =
\sum_{j=1}^{n} \varphi (\Delta t\, M)^{n-j}\, \dh_{j}.
\end{equation}
With the spatial truncation error
\[
\delta_h(t)=
\left(
  \begin{array}{c}
    \delta_h^L(t)\vspace{0.1cm} \\
    \delta_h^R(t) \\
  \end{array}
\right)
\]
as defined in Section \ref{conv}, one can express $b(t)$ as
\[
b(t) = u_h'(t)-Mu_h(t)-\delta_h(t).
\]
Inserting this into the definition of $\dh_j$ it readily
follows that
\begin{subeqnarray*}\label{}
\dh_j
&=& (I-\theta \Delta t\, M)^{-1} [(1-\theta)\,\Delta t\, \delta_h(t_{j-1})+
\theta\, \Delta t\, \delta_h(t_j) ] + \\
&&(I-\theta \Delta t\, M)^{-1} [u_h(t_j)-u_h(t_{j-1})-
(1-\theta)\,\Delta t\,u_h'(t_{j-1})-\theta\, \Delta t\,u_h'(t_j)].  % new
\end{subeqnarray*}
The above formula for $\dh_j$ consists of {\it two terms},
corresponding to the truncation error in space and the
truncation error in time.
In view of (\ref{epshat}), we shall study
$\varphi (\Delta t\, M)^{n-j}\, \dh_{j}$.\\

Concerning the first term, the bounds for the submatrices
of $\psi_{n-j}(\Delta t\, M) =
\varphi (\Delta t\, M)^{n-j}\,(I-\theta \Delta t\, M)^{-1} $
given by Lemma \ref{lemma psiM} directly lead to

\begin{eqnarray*}
&&|\psi_{n-j}(\Delta t\, M)\,
[(1-\theta)\,\Delta t\, \delta_h(t_{j-1})+
\theta\, \Delta t\, \delta_h(t_j) ]\,|_\infty
~\le\\\\
&&
\Delta t \cdot \left\{ K\cdot \max_{l=j-1,j} |\delta_h^L(t_l)|_\infty +
\frac{4(2K+1)S}{h_{m+2}}\cdot \max_{l=j-1,j} |\delta_h^R(t_l)|_\infty \right\}.
\end{eqnarray*}
Invoking the estimate (\ref{deltahR}) for $\delta_h^R$
this gives
\begin{eqnarray*}
&&|\psi_{n-j}(\Delta t\, M)\,
[(1-\theta)\,\Delta t\, \delta_h(t_{j-1})+
\theta\, \Delta t\, \delta_h(t_j) ]\,|_\infty
~\le\\\\
&&
\Delta t \cdot \left\{ K\cdot \max_{l=j-1,j} |\delta_h^L(t_l)|_\infty +
{\widehat \kappa} \cdot \max_{l=j-1,j} |\eta(t_l)|_\infty \right\},
\end{eqnarray*}
where ${\widehat \kappa} = (K+\tfrac{1}{2})\kappa$.\\
%with $\kappa$ and $\eta$ defined by (\ref{kappa_eta}).\\

Concerning the second term, for the function $g:[0,S] \rightarrow \R$
defined by
\[
g(s) = u(s,t_j)-u(s,t_{j-1})-
(1-\theta)\,\Delta t\,u_t(s,t_{j-1})-\theta\, \Delta t\,u_t(s,t_j)
%\quad (0\le s\le S)
\]
standard Taylor expansion shows that
% there exist constants $p$, $c_0$, $c_1$ such that
\[
|g(s)| \le c_0 (\Delta t)^{p+1}~~,~~|g'(s)| \le c_1 (\Delta t)^{p+1}
\quad (0\le s\le S)
\]
with
\begin{eqnarray*}
&& p=1~,~~c_0 = \tfrac{3}{2} ||u_{tt}||_\infty~,~~
~~c_1 = \tfrac{3}{2} ||u_{stt}||_\infty
\quad ~~({\rm if~}\tfrac{1}{2} < \theta \le 1), \\
&& p=2~,~~c_0 = \tfrac{5}{12} ||u_{ttt}||_\infty~,~~
c_1 = \tfrac{5}{12} ||u_{sttt}||_\infty
\quad ({\rm if~} \theta = \tfrac{1}{2}),
\end{eqnarray*}
where $||\cdot||_\infty$ designates the maximum norm of
a real function on $[0,S]\times [0,T]$.
Using now Lemma \ref{lemma psiMg} and the partitioning
\[
u_h(t)=
\left(
  \begin{array}{c}
    v_h(t)\vspace{0.1cm} \\
    w_h(t) \\
  \end{array}
\right)
\quad
{\rm with}~~v_h(t)\in\R^m~~{\rm and}~~w_h(t)\in\R^2,
\]
it follows that
\[
|\psi_{n-j}(\Delta t\, M)\, [u_h(t_j)-u_h(t_{j-1})-
(1-\theta)\,\Delta t\,u_h'(t_{j-1})-\theta\, \Delta t\,u_h'(t_j)]|_\infty  % new
\le c\, (\Delta t)^{p+1}
\]
with
\[
c=(4K+1)c_0 + (4K+2)Sc_1.
\]
\vskip0.2cm

Combining the above bounds, we are led to
\begin{eqnarray*}
|\eh_n|_\infty
&\le& \sum_{j=1}^{n} |\varphi (\Delta t\, M)^{n-j}\, \dh_{j}|_\infty\\
&\le& \sum_{j=1}^{n} \Delta t \cdot \left\{ K \cdot \max_{l=j-1,j} |\delta_h^L(t_l)|_\infty +
{\widehat \kappa} \cdot \max_{l=j-1,j} |\eta(t_l)|_\infty + c\, (\Delta t)^{p}\right\}\\
&\le& t_n \cdot \left\{ K \cdot\max_{0\leq\vartheta\leq t_n} |\delta_h^L(\vartheta)|_\infty
+{\widehat \kappa}\cdot\max_{0\leq\vartheta\leq t_n}\eta(\vartheta)+c\, (\Delta t)^{p}\right\}.
\end{eqnarray*}

\begin{flushright}
$\Box$
\end{flushright}\vskip0.2cm

%%%%%%%%%%%%%%%%%%%%%%%%%%%%%%%%%%%%%%%%%%%%%%%%%%%%%%%%%%%%%%%%%%%%%%%%%%%%%%%%%%%%%%
\section{Conclusions}\label{concl}
% new
In this paper we have analyzed the stability and convergence
of discretizations, both in space and time, of the
Black--Scholes PDE when it is provided with the linear
boundary condition.
% new
This condition
states that the second derivative of
the option value vanishes when the underlying asset price gets
large.
For the space discretization we considered finite difference
schemes and for the time discretization the well-known family
of $\theta$-methods.
% new
Concerning stability,
we derived
%, under natural assumptions,
tight inclusions for the maximum norm of $e^{tM}$ and
$\varphi (\Delta t\, M)^n$ where the matrix $M$ represents
the semidiscretized Black--Scholes PDE operator, $\varphi$
denotes the stability function of the $\theta$-method,
time $t\ge 0$, step size $\Delta t >0$ and integer $n\ge 0$.
The obtained inclusions reveal that these norms can grow
essentially directly proportional to the dimension of the
matrix $M$, i.e., the number of spatial grid points.
We subsequently proved the positive result that this growth
has in general no adverse effect on the convergence behavior
of the discretizations.

%The results in the present paper were derived for a particular
%discretization of the linear boundary condition.
In future research, in addition to the issues already
mentioned previously, we wish to extend the analysis to
the discretization of the linear boundary condition that
was considered by Windcliff, Forsyth \& Vetzal \cite{WFV04}.
% new
%This can be viewed as a variant of our discretization.
Also we aim at studying more advanced, multi-dimensional
PDEs in finance, such as the Heston PDE, as well as other
discretizations in space and time.

\section*{Acknowledgments}
The second author wishes to thank Tinne Haentjens for a
useful discussion on M-matrices.
She further acknowledges financial support by the Research
Foundation -- Flanders, FWO contract no.~1.1.161.10.N.

\end{document}